\documentclass[12pt]{article}
\usepackage[english]{babel}
\usepackage[T1]{fontenc}
\usepackage[latin1]{inputenc}
\usepackage{fullpage}
\usepackage{amssymb}
\usepackage{amsmath}
\usepackage{amsfonts}
\usepackage[all]{xy}
\usepackage{stmaryrd}

\begin{document}
\begin{center}
\LARGE{Positive Energy-Momentum Theorem for\\
AdS-Asymptotically Hyperbolic Manifolds}\\
\end{center}

\vspace{2pt}

\begin{center}
\large{\textsc{Daniel Maerten}}
\end{center}

\vspace{2pt}

\section{The Energy-Momentum}

\subsection{Introduction}

This paper proves a positive energy-momentum theorem under the (well known in general relativity) dominant energy condition, for AdS-asymptotically hyperbolic manifolds. An AdS-asymptotically hyperbolic manifold is by definition a manifold $(M,g,k)$ such that at infinity, the Riemannian metric $g$ and the symmetric 2-tensor $k$ tend respectively to the metric and second fundamental form of a standard hyperbolic slice of Anti-de Sitter (AdS).\\
Chru\'sciel and Nagy \cite{ChN} recently defined the notion of energy-momentum of an asymptotically hyperbolic manifold, which generalizes the analogous notion in the asymptotically flat case. Besides Chru\'sciel and Herzlich \cite{ChH} recently proved a positive mass theorem for asymptotically hyperbolic spin Riemannian manifolds (with zero extrinsic curvature).\\
The aim of the present paper is to extend this result to the non-zero extrinsic curvature case.

\subsection{Some Definitions and Notations}

We consider a Lorentzian manifold $N^{n+1}$ and a Riemannian spacelike hypersurface $M$. Using geodesic coordinates along $M$, we shall write  a neighbourhood of $M$ in $N$ as a subset of $]-\epsilon,\epsilon[ \times M$, endowed with the metric  $\gamma= -\text{d}t^{2} + g_{t}$. The Riemannian n-manifold $M$ has induced metric $g_{0}=g$ and second fundamental form $k:=(-\frac{1}{2}\frac{\text{d}}{\text{d}t}g_{t})_{|t=0}$. We assume that $(M,g,k)$ is AdS-asymptotically hyperbolic that is to say, the metric $g$ and the second fundamental form $k$ are asymptotic at infinity to the metric and the second fundamental form of a standard hyperbolic slice in AdS. More precisely we adopt the following\\

\noindent
\textbf{Definition.} $(M,g,k) $\textit{ is said to be AdS-asymptotically hyperbolic if there exists some compact K, a positive number R and a homeomorphism} $M\smallsetminus K \longrightarrow \mathbb{R}^{n}\smallsetminus B(0,R)$ \textit {called a chart at infinity such that in this chart we have}
$$\left\{ \begin{array}{lll}
e:=g-b=O(e^{-\tau r}), & \partial e= O(e^{-\tau r}), & \partial^{2} e = O(e^{-\tau r}),\\
k=O(e^{-\tau r}), & \partial k = O(e^{-\tau r}),
\end{array} \right.$$
\textit{for}  $\tau>n/2$ \textit{ and where} $\partial$ \textit{is taken with respect to the hyperbolic metric} $b= \text{d}r^{2} + \sinh^2 r g_{\Bbb S^{n-1}}$ \textit{with} $g_{\Bbb S^{n-1}}$ \textit{the standard metric of} $\Bbb S^{n-1}.$\\

\noindent
AdS space-time is merely denoted by $(\widetilde{N},\beta)$. If one considers  $\widetilde{N}$ as $\mathbb{R}^{n+1}$ then we will write $\beta=-\text{d}t^{2} +b_{t}$, with $b_{0}=b$ the hyperbolic metric.\\

%\noindent
%\textbf{Remark.} As we will see below, the definition of the energy-momentum only depends upon the behaviour of $(g,k)$ on  $M\smallsetminus K$. So we will assume that M and $\mathbb{H}^{3}$ are homeomorphic without any loss of generality. This also entails that  $N \cong \widetilde{N}$ as topological spaces. This asumption will be convenient in order to work on the whole manifold $N$ and not only on $]-\epsilon,\epsilon[ \times(M\smallsetminus K)$. \hspace{\stretch{4}}$\square$\\

\noindent
The motivation for the definition of the energy-momentum comes from the study of the constraints map which by definition is  
$$\begin{array}{lclc}
\Phi:& \mathcal{M}\times \Gamma(S^{2}T^{*}M)&\longrightarrow & C^{\infty}(M)\times \Gamma(T^{*}M)\\
& (h,p) & \longmapsto & \left(\begin{array}{c}\text{Scal}^{h} + (tr_{h}p)^{2}-\left|p\right|^{2}_{h}\\
2\left(\delta_{h}p + \text{d}tr_{h}p\right)
\end{array}\right),
\end{array}$$
where $ \mathcal{M}$ is the set of Riemannian metrics on the manifold $M$. Let us denote by $(\dot{h},\dot{p})$ an infinitesimal deformation of $(h,p)$. Now if we take a couple $(f, \alpha)\in C^{\infty}(M)\times \Gamma(T^{*}M)$ then we compute 
\begin{eqnarray*}
\left\langle (f, \alpha),(\Phi(h+\dot{h},p+\dot{p})- \Phi(h,p))\right\rangle &=& \delta ( f(\delta \dot{h} + \text{dtr}\dot{h})+i_{\nabla f}\dot{h}- (\text{tr}\dot{h})\text{d}f + 2i_{\alpha}\dot{p} - 2(\text{tr}\dot{p})\alpha)\\
&+& \delta ( <p,\dot{h}>\alpha +<h,\dot{h}>i_{\alpha}p -2i_{i_{\alpha}p}\dot{h})\\
&+& \left\langle \text{d}\Phi_{(h,p)}^{*}(f, \alpha),(\dot{h},\dot{p}) \right\rangle +Q(f, \alpha, h,p,\dot{h},\dot{k}),
\end{eqnarray*}
where <,> is the metric extended to all tensors, $\delta$ is the $h$-divergence operator, $\text{d}\Phi_{(h,p)}^{*}$ is the formal adjoint of the linearized constraints map at the point $(h,p)$, traces are taken with respect to $h$ and  $Q(f, \alpha,h,p,\dot{h},\dot{k})$ is a remainder which is linear with respect to $(f, \alpha)$ and at least quadratic with respect to $(\dot{h},\dot{p})$. Now considering the constraints map along the hyperbolic space embedded in AdS, that is to say $(h,k)=(b,0)$ and $(\dot{h},\dot{k})=(g-b=e,k)$ one finds
\begin{eqnarray*}
\left\langle (f, \alpha),(\Phi(g,k)- \Phi(b,0))\right\rangle &=& \delta ( f(\delta e + \text{dtr}e)+i_{\nabla f}e- (\text{tr}e)\text{d}f + 2i_{\alpha}k -2 (\text{tr}k)\alpha)\\
&+& \left\langle \text{d}\Phi_{(b,0)}^{*}(f, \alpha),(e,k) \right\rangle +Q(f, \alpha, b,k,e).
\end{eqnarray*}
As a consequence if we assume that $(M,g,k)$ is AdS-asymptotically hyperbolic and if the function $\left\langle (f, \alpha),(\Phi(g,k)- \Phi(b,0))\right\rangle$ is integrable on $M$ with respect to the measure $\text{dVol}_{b}$, then the energy-momentum $\mathcal{H}$ can be defined as a linear form on Ker $\text{d}\Phi_{(b,0)}^{*}$
$$\xymatrix{
 \mathcal{H}: (f,\alpha) \ar@{{|}->}[r] & \int_{S_{\infty}}- f(\delta_{b}e + \text{d}tr_{b}e) - i_{\nabla^{b} f}e +(tr_{b}e) \text{d}f - 2i_{\alpha^{\sharp}}k + 2(tr_{b}k)\alpha }.$$
The integrand in the formula of $\mathcal{H}$ is in index notation
$$ f( e_{j,i}^{i}-e_{i,j}^{i})-f^{,i}e_{ij}+(e_{i}^{i})f_{,j}- 2\alpha^{i}k_{ij}+2(k_{i}^{i})\alpha_{j},$$
where ``,'' stands for the $b$-derivatives and where $h^{i}=b^{ij}h_{j}$ for any tensor $h$. This integrand is the same as the one of Chru\'sciel and Nagy \cite{ChN} since each Killing vector fields on AdS is decomposable into the sum of some normal and  tangential components (with respect to a standard hyperbolic slice) which are in our case given by the couple $(f,\alpha)$ (see also \cite{HT}). More precisely, one can show, using Moncrief argument \cite{Mo}, that Ker $\text{d}\Phi_{(b,0)}^{*}\cong \mathfrak{Kill}(\text{AdS})$ where $ \mathfrak{Kill}(\text{AdS})$ denotes the Lie algebra of Killing vector fields on AdS,  since it satisfies the Einstein equations with a (negative) cosmological constant. The isometry group of AdS is O(n,2), and thereby $\mathfrak{Kill}(\text{AdS})\cong \mathfrak{so}(n,2)\cong N_{b}\oplus \mathfrak{so}(n,1)\cong N_{b}\oplus \mathfrak{Kill}(\Bbb H^{n}) $, where we have set $N_{b}=\{ f\in C^{\infty}(M) | \text{Hess}f=fb \}$. It is well known \cite{ChN}, \cite{ChH} that the application
$$\begin{array}{lcl}
\mathbb{R}^{n,1} & \longrightarrow & N_{b} \\
y_{k}            &\longmapsto & x_{k}:= y_{k|\mathbb{H}^{n}}
\end{array},$$
(where $(y_{k})_{k=0}^{n}$ are the standard coordinates) is an isometry, and the \textit{mass} part of the energy-momentum is a linear form on $N_{b}$ which is causal and positively oriented as soon as $\text{Scal}^{g} \geq -n(n-1)=\text{Scal}^{b}$. Remark that the sharpest integrability conditions in order to make $\mathcal{H}$ well defined and invariant under asymptotic isometries have been found by Chru\'sciel and Nagy still in \cite{ChN}. However for the sake of simplicity one can use instead of the integrability condition  $\left\langle (f, \alpha),(\Phi(g,k)- \Phi(b,0))\right\rangle \in L^{1}(M, \text{dVol}_{b})$, the less general but more convenient condition $\left|\Phi(g,k)-\Phi(b,0) \right|e^{r} \in  L^{1}(M, \text{dVol}_{b})$.\\

\noindent
\textbf{Remark.} In the asymptotically flat situation, the energy-momentum is also a linear form on $\mathbb{R}^{n,1}\oplus \mathfrak{so}(n,1)$ where the first component corresponds to {\it translational} isometries and the second one to {\it rotations}. This interpretation gave rise to the respective terminology of linear and angular momentum. In the AdS-asymptotically hyperbolic situation, one cannot identify some linear momentum in the decomposition $\mathfrak{so}(n,2)\cong \mathbb{R}^{n,1}\oplus \mathfrak{so}(n,1)$, since the  first component $\mathbb{R}^{n,1}$ in $\mathfrak{so}(n,2) $ is not of translational nature. This whole first component of the energy is then called the {\it mass functional}  and it only remains some {\it angular momentum}. Moreover physicists often call the limit of integrals $\mathcal{H}(f,\alpha )$ {\it global charges} and so the positive energy-momentum theorem could be consequently renamed global inequalities theorem. Some supplementary details on the physical interpretation of our result can be found in the forthcoming note by Chru\'sciel and the author \cite {ChM}. \hspace{\stretch{4}}$\square$\\

\subsection{Statement of the Theorems and Comments}

As a matter fo fact, we know that, given a chart at infinity, $\mathcal{H}$ can be considered as a vector of $\mathbb{R}^{n,1}\oplus \mathfrak{so}(n,1)$ and will be denoted by $M\oplus\Xi$. The vector $M$ is the mass part \cite{ChH} of $\mathcal{H}$, and $\Xi$ is the angular momentum. We will prove the existence of  a Hermitian application 
$$Q :\xymatrix{
\Bbb C^{d} \ar[r]^-{\mathcal{K}} & \Bbb R^{n,1}\oplus \mathfrak{so}(n,1)  \ar[r]^-{\mathcal{H}} & \Bbb R
} ,$$
which has to be non-negative under some energy assumptions but nonetheless $Q$ quite difficult to explicite in general.\\
However, in dimension $n=3$, we can be more specific giving the explicite formula of $Q$ in terms of the components $M\in\mathfrak{M}\subset \text{M}(2,\Bbb C)$ (cf. section 2.4 for the definition of $\frak M$) and $\Xi\in \mathfrak{sl}(2,\Bbb C)$ of the energy-momentum $\mathcal{H}$. More precisely we will show that 
$$Q=2
\left(\begin{array}{cc}
	\widehat{M}& \Xi\\
	\Xi^{*} & M
\end{array}\right),$$
where  $\widehat{M}$ is the transposed comatrix of $M$. We will also treat the case where the slice $M$ has a compact inner boundary $\partial M$ whose induced metric and second fundamental form are respectively denoted by  $\breve{g}$ and $\breve{k}$. To this end,  we have to define the vector field  $\vec{k}:=(-\text{tr}\breve{k}+(n-1))e_{0}+k(\nu)$ along the boundary $\partial M$. We can now state the\\

\noindent
\textbf{Positive Energy-Momentum Theorem.} \textit{Let} $(M^{n},g,k)$ \textit{be an AdS-asymptotically hyperbolic spin Riemannian manifold satisfying the decay conditions stated in section 1.2 and the following conditions\\
(i)} $\left\langle (f, \alpha),(\Phi(g,k)- \Phi(b,0))\right\rangle \in L^{1}(M,\text{dVol}_{b})$ \textit{for every} $(f,\alpha)\in N_{b}\oplus \mathfrak{Kill}(M,b),$ \\
\textit{(ii) the relative version of the dominant energy condition (cf. section 2.2) holds, that is to say} $(\Phi(g,k)- \Phi(b,0))$ \textit{is a positively oriented causal (n+1)-vector along $M$,\\
(iii) in the case where M has a compact boundary $\partial M$, we assume moreover that $\vec{k}$ is causal and  positively oriented along $\partial M$.\\
Then there exists  a (hardly explicitable) map} $\Bbb R^{n,1} \oplus \mathfrak{so}(n,1) \longrightarrow  \text{Herm}(C^d)$ {\it which sends, under the assumptions (i-iii), the energy-momentum on a non-negative Hermitian form Q.\\
Moreover, when n=3, we can explicite Q in terms of the components of the energy-momentum as described above.}\\

\noindent
Classical algebra results give the non negativity of each principal minors of $Q$ which provide a set of inequalities on the coefficients of $\mathcal{H}$ that are explicitely written in the appendix in dimension $n=3$ (cf. section 5.1).\\
This result is new (even though many formal arguments where given by Gibbons, Hull and Warner in \cite{GHW}) and  based on the recent {\it global charge} definition of Chru\'sciel and Nagy for AdS-asymptotically hyperbolic manifolds, which comes from the Hamiltonian description of General Relativity. Our approach is purely Riemannian, the Lorentzian connection and manifold introduced are auxiliary since everything is restricted to the Riemannian slice $M$  (contraryly to \cite{GHW}). In the other hand, the positive mass theorem for Minkowski-asymptotically hyperbolic initial data sets of Chru\'sciel, Jezierski and ~\L\c{e}ski \cite{ChJL} is also different from ours, since their Riemannian hypersurface is supposed to be asymptotic at infinity to a standard hyperbolic slice of Minkowski space-time (in that case the extrinsic curvature does not tend to 0) . Then considering the translational Killing vector fields of Minkowski, they defined a hyperbolic 4-momentum  (usually called {\it Trautman-Bondi mass}) and proved that it is timelike and future directed under the dominant energy condition (and some other technical assumptions).\\
Remark finally that our result extends the positive mass theorem of Chru\'sciel and Herzlich \cite{ChH} in dimension $n$, since if one supposes that $k=0$ then we recover their result: the mass functional $M$ has to be time-like future directed.\\
As regards the rigidity part we have the\\

\noindent
\textbf{Theorem.} {\it Under the assumptions of the positive energy-momentum theorem},  ${\text{tr}Q=0}$ {\it implies that  $(M,g,k)$ is  isometrically embeddable in} AdS$^{n,1}$.\\

\noindent
This result is optimal in the sense that one could not reasonably hope better than being able to embbed isometrically our triple $(M,g,k)$ in AdS. \\
Some additional but partial results will be proved (also for the Trautman-Bondi 4-momentum) in order to weaken the defining condition of rigidity.

\subsection{Organisation of the Paper}

In section 2, we give the necessary geometric background by recalling some basic facts on spinors and defining the Killing connection used in the remainder of the paper. We also prove the Bochner-Lichnerowicz-Weitzenböck-Witten formula with respect to our Killing connection and deduce an integration formula.\\
In section 3, we prove the positive energy-momentum theorem: we remark that the boundary contribution of the integrated Bochner-Lichnerowicz-Weitzenböck-Witten formula can be identified to the {\it global charges} $\mathcal{H}(f,\alpha )$, for some choices of $(f,\alpha )$. This can be done using the same ideas as in \cite{ChH} but in a Lorentzian situation, and extends the quite technical computations of  \cite{ChH} in a non-trivial way, since the algebraic structures are different (spinors, Hermitian scalar product, gauge etc...) and since new terms (involving the extrinsic curvature) appeared and had to be identified. We then make the analysis of the Dirac operator (we also treat the case where $M$ has a compact boundary) which gives the non-negativity of the {\it globlal charges}  $\mathcal{H}(f,\alpha )$ when the couple $(f,\alpha )$ comes from an imaginary Killing spinor of AdS$^{n,1}$. Then we restrict to dimension $n=3$, and completely study the imaginary Killing spinors of AdS$^{3,1}$ in order to interpret the non-negativity  of the {\it global charges} as the non-negativity of the Hermitian matrix $Q$ on $\Bbb C^{4}$.\\
Section 4 is devoted to the proof of the {\it rigidity results}.\\
The last section is an appendix which gives the non-negativity of $Q$ in dimension $n=3$ seen through its coefficients, and proves some {\it rigidity results} for the Trautman-Bondi mass \cite{ChJL}.

\section{Geometric Background}

All the definitions and conventions of this section will be for any $n\geq 3$, where $n$ is the dimension of the AdS-asymptotically hyperbolic slice, except if the dimension is explicitely mentioned to be 3. 

\subsection{Connections and curvatures}

$\nabla, \overline{\nabla}$ denote respectively the Levi-Civita connections of $\gamma$ and \textit{g}. Let us take a spinor field $\psi \in \Gamma(\Sigma)$ and a vector field $X\in\Gamma(TM)$, then 
$$\left\{\begin{array}{cll} 
\nabla_{X}\psi & = & \overline{\nabla}_{X}\psi - \frac{1}{2}k(X)\cdot e_{0} \cdot \psi \\
\left\langle k(X),Y \right\rangle_{\gamma} &= & \left\langle \nabla_{X}Y,e_{0}\right\rangle_{\gamma}
\end{array} \right. .$$
In these formulae $\cdot$ denotes the Clifford action with respect to the metric $\gamma$, and $e_{0}=\partial_{t}$. We will use different notations when we have to make the difference between the Clifford action with respect to the metric $\gamma$ or $\beta$.\\

\noindent
\textbf{Definition.} {\it The Killing equation on a spinor field $\tau\in \Gamma(\Sigma)$ is 
$$ \widehat{D}_{X}\tau:=D_{X} \tau + \frac{i}{2}X \cdot_{\beta } \tau =0  \quad \forall X \in \Gamma (T M),$$
where $D$ denotes the Levi-Civita connection of AdS along $M$. Such a $\widehat{D}$-parallel spinor field is called a $\beta $-imaginary Killing spinor and we denote $\tau \in IKS(\Sigma )$. In the same way, a  $\widehat{\nabla}$-parallel spinor field (where $ \widehat{\nabla}_{X} := \nabla_{X} +  \frac{\textbf{i}}{2} X\cdot_{\gamma }$) is called a $\gamma $-imaginary Killing spinor.}\\
 
\noindent
Notice that the equation $\widehat{D}\tau=0$ is neither the Killing equation in AdS nor in $\Bbb{H}^{n}$, but the Killing equation in AdS along $\mathbb{H}^{n}$ (in particular the imaginary Killing spinors considered here are not the one of \cite{AD}, \cite{ChH}).\\
Now  if  $R, \widehat{R}$ are the respective curvatures of $\nabla$ and $\widehat{\nabla}$, we have the relation
$$ \widehat{R}_{X,Y} = R_{X,Y} -\frac{1}{4} (X\cdot Y-Y\cdot X)\cdot  \quad ,$$
where we use the convention of \cite{LMi} for the curvature.

\subsection{Bochner-Lichnerowicz-Weitzenböck-Witten Formula and the Dominant Energy Condition}

From now on $\left(e_{k}\right) _{k=0}^{n}$ is an orthonormal basis at the point with respect to the metric $\gamma$. We define the Dirac-Witten operators

$$\mathfrak{D}\psi=\sum^{n}_{k=1}e_{k}\cdot\nabla_{e_{k}}\psi, \quad \widehat{\mathfrak{D}}\psi=\sum^{n}_{k=1}e_{k}\cdot\widehat{\nabla}_{e_{k}}\psi,$$
where $n$ is the dimension of the spacelike slice.\\

\noindent
\textbf{Lemma.}\textit{(Bochner-Lichnerowicz-Weitzenböck-Witten formula)}
$$\widehat{\mathfrak{D}}^{*}\widehat{\mathfrak{D}}
=\widehat{\nabla}^{*}\widehat{\nabla} +\widehat{\mathfrak{R}},$$
\textit{where} $ \widehat{\mathfrak{R}}:=\frac{1}{4} \left(\text{Scal}^{\gamma} + n(n-1) + 4 \text{Ric}^{\gamma}(e_{0},e_{0}) + 2e_{0}\cdot \text{Ric}^{\gamma}(e_{0}) \right).$\\

\noindent
\textsc{Proof.} The Dirac-Witten operator $\mathfrak{D}$ is clearly formally self adjoint, and we have the classical Bochner-Lichnerowicz-Weitzenböck formula (cf. \cite{H},\cite{PT} for instance) ${\mathfrak{D}^{*}\mathfrak{D}=\mathfrak{D}^{2}= \nabla^{*}\nabla + \mathfrak{R}}$,
where  $ \mathfrak{R}:=\frac{1}{4} \left(\text{Scal}^{\gamma} + 4 \text{Ric}^{\gamma}(e_{0},e_{0}) + 2e_{0}\cdot \text{Ric}^{\gamma}(e_{0}) \right)$. We also know that $\widehat{\mathfrak{D}}=\mathfrak{D}-\textbf{i}\frac{n}{2}$ and so we get 
$$\widehat{\mathfrak{D}}^{*}\widehat{\mathfrak{D}}=\nabla^{*}\nabla + \mathfrak{R}+ \frac{n^{2}}{4}, $$
but finally remarking that $ \widehat{\nabla}^{*}\widehat{\nabla}= \nabla^{*}\nabla + \frac{n}{4} $, we obtain our formula.\hspace{\stretch{4}}$\square$\\
We derive an integration formula from the Bochner-Lichnerowicz-Weitzenböck-Witten identity considering the 1-form $\theta$ on $M$ defined by $\theta(X)= \left\langle \widehat{\nabla}_{X}\psi + X\cdot \widehat{\mathfrak{D}}\psi,\psi   \right\rangle_{\gamma}$, where $\psi$ is a spinor field. Straightforward computations lead to the following $g$-divergence formula
$$ \text{div}\theta =  \left\langle\widehat{\mathfrak{D}}\psi ,\widehat{\mathfrak{D}}\psi \right\rangle_{\gamma}- \left\langle\widehat{\mathfrak{R}}\psi,\psi \right\rangle_{\gamma} -  \left\langle \widehat{\nabla}\psi,\widehat{\nabla}\psi \right\rangle_{\gamma}.$$
Let $S_{r}$ the $g$-geodesic sphere of radius $r$ and centered in a point of M. The radius $r$ is supposed to be as large as necessary. We denote by $M_{r}$ the interior domain of $S_{r}$ and $\nu_{r}$ the (pointing outside) unit normal. Integrating our divergence formula over $M_{r}$ and using Stokes theorem, we get
$$\int_{M_{r}} \left|\widehat{\mathfrak{D}}\psi\right|^{2}_{\gamma} = 
\int_{M_{r}}\left(\left|\widehat{\nabla}\psi\right|^{2}_{\gamma}+\left\langle \widehat{\mathfrak{R}}\psi,\psi\right\rangle_{\gamma}\right)
-\int_{S_{r}}\left\langle \widehat{\nabla}_{\nu_{r}}\psi+ \nu_{r}\cdot\widehat{\mathfrak{D}}\psi, \psi\right\rangle_{\gamma} \text{dVol}_{S_{r}}  \quad .$$

\noindent
Let us now consider the Einstein tensor $G=\text{Ric}^{\gamma}-\frac{1}{2}\text{Scal}^{\gamma}\gamma$ with respect to the metric $\gamma$. The dominant energy condition \cite{Wald} says that the speed of energy flow of matter is always less than the speed of light. More precisely, for every positively oriented time-like vector field $v$, the energy-momentum current of density of matter $-G(v,.)^{\sharp}$ must be time-like or null, with the same orientation as $v$. The assumption we make in order to prove the positive energy-momentum theorem is a relative version of the dominant energy condition: {\it $ -\left(G-\frac{n(n-1)}{2}\gamma\right)(e_{0})$ is a positively oriented time-like or null vector along $M$}. Some easy computations give
\begin{eqnarray*}
\text{Scal}^{\gamma} & = & 2(G(e_{0},e_{0})-\text{Ric}^{\gamma}(e_{0},e_{0})) \\
 e_{0}\cdot \text{Ric}^{\gamma}(e_{0})&=& e_{0}\cdot G_{|TM}(e_{0})-\text{Ric}^{\gamma}(e_{0},e_{0}),
\end{eqnarray*}
where $G_{|TM}(e_{0})= \sum_{k=1}^{3}G(e_{0},e_{k})e_{k}$. Thereby
\begin{eqnarray*}
	\widehat{\mathfrak{R}}&=& \frac{1}{4}\left( 2G(e_{0},e_{0}) +\left(n(n-1)\right)+2 e_{0}\cdot G_{|TM}(e_{0})  \right)\cdot\\
	&=& \frac{1}{2}\left(\left(  G(e_{0},e_{0})+\frac{n(n-1)}{2} \right)e_{0} -  G_{|TM}(e_{0})  \right) \cdot e_{0}\cdot\\
	&=& \frac{1}{2} \left( \left(G(e_{0},e_{0})-\frac{n(n-1)}{2}\gamma(e_{0},e_{0})\right)e_{0} - G_{|TM}(e_{0})  \right)\cdot e_{0}\cdot \\
	&=& - \frac{1}{2}\left(G-\frac{n(n-1)}{2}\gamma\right)(e_{0})\cdot e_{0}\cdot \quad .
\end{eqnarray*}
Our assumption gives the non negativity of the spinorial endomorphism $\widehat{\mathfrak{R}}$ that is to say $ \left\langle \widehat{\mathfrak{R}}\psi,\psi \right\rangle \geq 0$ for every spinor field $\psi$.\\

\noindent
\textbf{Remark.} We can express the dominant energy condition in terms of the constraints as in section 1.3  since
$$ -\left(G-\frac{n(n-1)}{2}\gamma\right)(e_{0})=\frac{1}{2}\left(\Phi(g,k)-\Phi(b,0)\right).$$\hspace{\stretch{4}}$\square$\\

\subsection{Spinorial Gauge}

In the same way as Andersson and Dahl \cite{AD}, but in a Lorentzian situation, we compare spinors in $\Sigma$ (along $M$) with respect to the two different metrics $\beta$ and $\gamma$. This can be done according to \cite{BG} as soon as the tubular neighbourhood of $M$ in $N$ is small enough. Consequently we suppose that both metrics are written in Gaussian coordinates $ \beta= - \text{d}t^{2} + g_{t}, \ \gamma=- \text{d}t^{2} + b_{t}$ on $]-\epsilon ,+\epsilon [ \times M$ for $\epsilon $ small enough. We define the spinorial gauge $\mathcal{A}\in \Gamma (\text{End}(\mathbb{T}))$ with the relations
$$\left\{\begin{array}{lll}
\gamma(\mathcal{A}X,\mathcal{A}Y)&=& \beta(X,Y)\\
\gamma(\mathcal{A}X,Y)&=& \gamma(X,\mathcal{A}Y)	
\end{array}\right. ,$$
where $\mathbb{T}$ is $TN$ restricted to $M$. The first relation says that $\mathcal{A}$ sends $\beta$-orthonormal frames on $\gamma$-orthonormal frames whereas the second one means that the endomorphism $\mathcal{A}$ is symmetric. We notice that these relations are only satisfied along $M=\left\{ t=0 \right\}$ and can also be written in the following way
$$\left\{\begin{array}{rll}
\mathcal{A}e_{0} &=& e_{0} \\
g(\mathcal{A}X,\mathcal{A}Y)&=& b(X,Y)\\
g(\mathcal{A}X,Y)&=& g(X,\mathcal{A}Y)	
\end{array}\right. .$$
Consequently $\mathcal{A}$ is an application $\text{P}_{\text{SO}_{0}(n,1)}(\beta)_{|M} \longrightarrow \text{P}_{\text{SO}_{0}(n,1)}(\gamma)_{|M}$, which can be covered by an application still denoted $\mathcal{A}:\text{P}_{\text{Spin}_{0}(n,1)}(\beta)_{|M} \longrightarrow \text{P}_{\text{Spin}_{0}(n,1)}(\gamma)_{|M}$. This application carries $\beta$-spinors on $\gamma$-spinors so that we have the compatibility relation about the Clifford actions of $\beta$ and $\gamma$ 
$$\mathcal{A}(X \cdot_{\beta}\sigma)= (\mathcal{A}X)\cdot_{\gamma}(\mathcal{A}\sigma),$$
for every $X \in \Gamma (\mathbb{T}), \sigma \in \Gamma(\Sigma)$ and where $\cdot_{\beta}, \cdot_{\gamma}$ denotes the Clifford actions respectively of $\beta$ and $\gamma$. Remark that our gauge is more sophisticated that the one of \cite{AD} since it deals with the trace of Lorentzian structures (metrics, spinors, Hermitian scalar product etc...) along the spacelike slice $M$.\\
We define a new connection $\widetilde{\nabla}X= \mathcal{A}(\overline{D}\mathcal{A}^{-1}X)$ along $M$. It is easy to check that $\widetilde{\nabla}$ is $g$-metric and has torsion $\widetilde{T}(X,Y)=-((\overline{D}_{X}\mathcal{A})\mathcal{A}^{-1}Y-(\overline{D}_{Y}\mathcal{A})\mathcal{A}^{-1}X)$.
We extract some formulae for later use
$$2 g\left( \widetilde{\nabla}_{X}Y-\overline{\nabla}_{X}Y,Z \right)= 
g\left(\widetilde{T}(X,Y),Z \right) -g\left(\widetilde{T}(X,Z),Y \right) -g\left(\widetilde{T}(Y,Z),X \right).$$
Now we intend to compare the connexions $\overline{\nabla}$ and $\widetilde{\nabla}$ on $\Sigma$. $(\sigma_{s})_{s}$ denotes the spinorial frame corresponding to the orthonormal frame $(e_{k})^{n}_{k=0}$, and $\overline{\omega},\widetilde{\omega}$ are the connection 1-forms respectively of $\overline{\nabla}$ and $\widetilde{\nabla}$
$$\begin{array}{lll}
\overline{\omega}_{ij}&=& g(\overline{\nabla} e_{i}, e_{j})\\
\widetilde{\omega}_{ij}&=& g(\widetilde{\nabla} e_{i}, e_{j}),
\end{array}$$
and if we take a general spinor $\varphi= \varphi^{s}\sigma_{s}$, their derivatives are given by
\begin{eqnarray*}
	\overline{\nabla}\varphi = \text{d}\varphi^{s}\otimes\sigma_{s} + \frac{1}{2}\sum_{i<j}\overline{\omega}_{ij}\otimes e_{i}\cdot_{\gamma}e_{j}\cdot_{\gamma}\varphi \\
	\widetilde{\nabla}\varphi = \text{d}\varphi^{s}\otimes\sigma_{s} + \frac{1}{2}\sum_{i<j}\widetilde{\omega}_{ij}\otimes e_{i}\cdot_{\gamma}e_{j}\cdot_{\gamma}\varphi,
\end{eqnarray*}
and as a consequence
$$ (\overline{\nabla} -\widetilde{\nabla})\varphi = \frac{1}{4}\sum_{i,j=0}^{n}(\overline{\omega}_{ij}-\widetilde{\omega}_{ij}) \otimes e_{i}\cdot_{\gamma}e_{j}\cdot_{\gamma}\varphi  \quad .$$

\subsection{Tangent and Spinor Bundles}

In this paper, the model spaces AdS$^{n,1}$ and $\mathbb{H}^{n}$ are considered as symmetric spaces:
$$\xymatrix{\mathbb{H}^{n}=\text{Spin}_{0}(n,1)/\text{Spin}_{n} \ar@{^{(}->}[r]& \text{AdS}^{n,1}=\text{Spin}_{0}(n,2)/\text{Spin}_{0}(n,1)},$$
so that every section of any natural fiber bundle above $\mathbb{H}^{n}$ can be seen as a function  $\text{Spin}_{0}(n,1)\longrightarrow \Bbb C^{d}$ which is $\text{Spin}(n)$--equivariant (with $d$ depending upon $n$). We can be more explicite when we take $n=3$ (this fact is due to the exceptional isomorphisms of Lie groups below)
$$\xymatrix{\mathbb{H}^{3}=\text{SL}(2,\mathbb{C})/\text{SU}(2) \ar@{^{(}->}[r]& \text{AdS}^{3,1}=\text{Spin}_{0}(3,2)/\text{Spin}_{0}(3,1)},$$
with $\text{SU}(2) \cong \text{Spin}(3) $ and $\text{SL}(2,\mathbb{C})\cong \text{Spin}_{0}(3,1) $.\\
The spinor bundle of AdS is  $\Sigma_{AdS} = \text{Spin}_{0}(3,2) \times _{\tilde{\rho}} \mathbb{C}^{4}$, where $\text{Spin}_{0}(3,2)$ is the bundle of the $\text{Spin}_{0}(3,1)$-frames in AdS, and $\tilde{\rho}$ is the standard  représentation of $\text{SL}(2,\mathbb{C})$ on $\mathbb{C}^{4}\cong \mathbb{C}^{2} \oplus \overline{\mathbb{C}^{2}}'$. In other words
$$
\begin{array}{cccl}
\tilde{\rho} : & \text{SL}(2,\mathbb{C})& \longrightarrow & \text{M}_{4}(\mathbb{C})\\
 & \tilde{g} & \longmapsto & 
 \left(
\begin{array}{cc}
\tilde{g}& 0\\
0& (\tilde{g}^{*})^{-1}	
\end{array}\right)
	\end{array},$$
where $A^{*}=^{t}\overline{A},  A \in \text{M}_{2}(\mathbb{C})$.
When we restrict this bundle to the hypersurface $\mathbb{H}^{3}$ we have $\Sigma = \text{SL}(2,\mathbb{C}) \times _{\tilde{\rho}_{|\text{SU}(2)}} \mathbb{C}^{4}$.\\

\noindent
\textbf{Proposition.} $\Sigma$ \textit{and} $\mathbb{H}^{3}\times \mathbb{C}^{4}$ \textit{are isomorphic thanks to the following trivialisation:}
$$
\begin{array}{cccl}
T: & 	\Sigma & \longrightarrow & \mathbb{H}^{3}\times \mathbb{C}^{4}\\
 & \{\tilde{e},w\} & \longmapsto & \left([\tilde{e}],\tilde{\rho}(\tilde{e})w \right)
\end{array},$$
\textit{where} $\{\tilde{e},w\}$ \textit{denotes the class of} $(\tilde{e},w)\in \text{SL}(2,\mathbb{C})\times \mathbb{C}^{4}$ \textit{in} $\Sigma$ \textit{, and} $[\tilde{e}]$ \textit{denotes the class of} $\tilde{e} \in \text{SL}(2,\mathbb{C})$ \textit{in} $\mathbb{H}^{3}=\text{SL}(2,\mathbb{C})/\text{SU}(2)$.\\

\noindent
The construction of $\mathbb{T}_{AdS}$, the tangent bundle of AdS, is quite similar to the construction of the spinor bundle. Still noticing that the principal bundle of $\text{SO}_{0}(3,1)$-frames in AdS is isomorphic to $\text{SO}_{0}(3,2)$, we write $\mathbb{T}_{AdS} =\text{SO}_{0}(3,2) \times_{\rho}\mathbb{R}^{4}$, where $\rho$ is the standard representation of $\text{SO}_{0}(3,1)$ on $\mathbb{R}^{4}$. By restriction to the hypersurface $\mathbb{H}^{3}$, we obtain $ \mathbb{T}=\text{SO}_{0}(3,1) \times_{\rho_{|\text{SO}(3)}}\mathbb{R}^{4}$, where $\text{SO}(3)$ is by definition the isotropy group of $f_{0}$ if $\left(f_{k}\right)^{3}_{k=0}$ denotes the canonical basis of $\mathbb{R}^{4}$.\\

\noindent
\textbf{Proposition.} $\mathbb{T}$ \textit{and} $\mathbb{H}^{3}\times \mathbb{R}^{4}$ \textit{are isomorphic thanks to the following trivialisation:}
$$
\begin{array}{cccl}
T: & 	\mathbb{T} & \longrightarrow & \mathbb{H}^{3}\times \mathbb{R}^{4}\\
 & \{e,u\} & \longmapsto & \left([e],\rho(e)u \right)
\end{array},$$
\textit{where }$\{e,u\}$ \textit{denotes the class of} $(e,u)\in \text{SO}_{0}(3,1)\times \mathbb{R}^{4}$ \textit{in} $\mathbb{T}$\textit{, and} $[e]$ \textit{denotes the class of} $e \in \text{SO}_{0}(3,1)$ \textit{in} $\mathbb{H}^{3}=\text{SO}_{0}(3,1)/\text{SO}(3)$.\\

\noindent
We are going to define the Clifford action on $\Sigma$, in the same way as in \cite{PT}. To this end, we denote by  $(\mathbb{R}^{4},q)$ the Minkowski space-time of signature (3,1), where $q= -\text{d}y^{2}_{0}+\text{d}y^{2}_{1}+\text{d}y^{2}_{2}+\text{d}y^{2}_{3}$. This space  is isometric to a subspace of $\text{M}_{2}(\mathbb{C})$ via
$$\begin{array}{lclc}
\Lambda: &   (\mathbb{R}^{4},q)& \longrightarrow & 	\mathfrak{M}:=\left(\left\{A \in \text{M}_{2}(\mathbb{C})/ A^{*}=A  \right\}, -\text{det}\right)\\
&  y=\left(y_{i}\right)^{3}_{i=0} & \longmapsto & \left(
\begin{array}{ll}
y_{0}+y_{1} & y_{2}+iy_{3} \\
y_{2}-iy_{3} & y_{0}-y_{1}	
\end{array}\right)
\end{array}.$$ 
We have thus the following real vector space isomorphisms:
$$\begin{array}{lll}
\text{M}_{2}(\mathbb{C}) & \cong &  \mathfrak{u}(2)\oplus \mathfrak{M}\\
\mathfrak{sl}_{2}(\mathbb{C}) & \cong &  \mathfrak{su}(2) \oplus \left( \mathfrak{M}\cap \mathfrak{sl}_{2}(\mathbb{C})\right)\\
 & \cong & \mathfrak{su}(2) \oplus \mathfrak{G},
\end{array}$$
and $\mathfrak{G}\cong \mathbb{R}^{3}$.  In order to make the value of the sectional curvature of $\mathbb{H}^{3}$ equal to -1, when we consider $\mathbb{H}^{3}=\text{SL}(2,\mathbb{C})/\text{SU}(2)$ as a symmetric space, we have to consider $\mathbb{R}^{4}$ endowed with $4q$ and not $q$, and consequently the embedding of the Clifford algebra $\text{C}\ell_{3,1}$ in $\text{M}_{4}(\mathbb{C})$ becomes
$$ \Theta: X \in  \mathfrak{M} \longmapsto \left(\begin{array}{cc}
0 & 2X \\
2\widehat{X} & 0	
\end{array}\right),$$
where $\widehat{X}$ means the transposed comatrix of X.\\
It will be convenient to see $\mathbb{T}$ as $\text{SL}(2,\mathbb{C})\times_{\mu} \mathfrak{M}$, where $\mu$ is the universal covering of $\text{SO}_{0}(3,1)$ by $\text{SL}(2,\mathbb{C})$, which is given by:
$$\begin{array}{cccc}
\mu :  & 	\text{SL}(2,\mathbb{C}) & \longrightarrow & \text{SO}_{0}(3,1)\\
 & \tilde{g} &\longmapsto & \left(\tilde{g}: X \in\mathfrak{M}\mapsto \tilde{g}X\tilde{g}^{*}\right)
\end{array}.$$
We can now define the Clifford action. Let us take $e \in \text{SO}_{0}(3,1)$ and $\tilde{e} \in \text{SL}(2,\mathbb{C})$ such that $e=\mu(\widetilde{e})$. A vector $X=X[e]$ tangent at the point $[e]=[\tilde{e}] \in \mathbb{H}^{3}$, is a class $\{e,u\} \in \mathbb{T}$. A spinor $\sigma = \sigma [\tilde{e}]$ at the same point is likewise a class $\{\tilde{e},w\} \in \Sigma$. The result of the Clifford action of $X$ on $\sigma$ is the spinor $(X \cdot \sigma)[\tilde{e}] = \{e,u\} \cdot \{\tilde{e},w\} = \{\tilde{e},\Theta(u)w\}$. We define a sesquilinear inner product (not definite positive) $(\cdot,\cdot)$ on $\mathbb{C}^{4}\cong \mathbb{C}^{2} \oplus \overline{\mathbb{C}^{2}}'$ as in \cite{PT} $(\xi,\eta):= \left\langle \xi_{1},\eta_{2} \right\rangle_{\mathbb{C}^{2}} +\left\langle \xi_{2},\eta_{1} \right\rangle_{\mathbb{C}^{2}}$, where $ \xi= \binom{\xi_{1}}{\xi_{2}}, \eta=\binom{\eta_{1}}{\eta_{2}} \in \mathbb{C}^{4}$ and where $\left\langle\cdot ,\cdot\right\rangle_{\mathbb{C}^{2}}$  is the standard Hermitian product on $\mathbb{C}^{2}$. This induces a sesquilinear product on $\Sigma$ by $ \left(\{\tilde{e}, \xi \},  \{\tilde{e}, \eta \}\right):= (\xi,\eta)$. In the same way we define a scalar product on $\Sigma$ setting
$$\begin{array}{lcl}
\left\langle \{\tilde{e}, \xi \},  \{\tilde{e}, \eta \}\right\rangle &:=& \left(\frac{1}{2}f_{0}\cdot\{\tilde{e}, \xi \},  \{\tilde{e}, \eta \}\right)\\
 & = & (\{\tilde{e}, \frac{1}{2}\Theta(f_{0})\xi \},  \{\tilde{e}, \eta \})\\
 & = & \left\langle \xi,\eta \right\rangle_{\mathbb{C}^{4}},
\end{array}$$
where $\left\langle \cdot,\cdot\right\rangle_{\mathbb{C}^{4}}$ denotes the standard Hermitian product on $\mathbb{C}^{4}$.\\
Since $\text{SL}(2,\mathbb{C})$ is the 2-sheeted covering of $\text{SO}_{0}(3,1)$, there exists a natural (left) action of $\text{SL}(2,\mathbb{C})$ on $\Sigma$ which is derived from the natural (left) action of $\text{SO}_{0}(3,1)$ on $\mathbb{T}$: the action of the group of the isometries of AdS preserving the slice $\mathbb{H}^{3}$ that is $\tilde{g} \ast \{\tilde{e},w\}= \{\tilde{g}\tilde{e}, w\}$,
with $\tilde{g} \in \text{SL}(2,\mathbb{C})$ and  $\sigma[\tilde{e}]=\{\tilde{e},w\}$ a spinor at $[\tilde{e}]$. To have the action on a section  $\sigma \in\Gamma(\Sigma)$ we set as usual $(\tilde{g} \ast \sigma )[\tilde{e}]= \tilde{g}*\sigma (\tilde{g}^{-1}\tilde{e})$.\\

\section{Positive Energy-Momentum Theorem}

In this section, the dimension will be $n\geq3$ expect if $n$ is explicitely mentioned to be 3.\\
Moreover $f$ will denote a smooth cutoff function which is 0 on $M$ except on a small neighbourhood of the infinity boundary of $M$ where $f\equiv 1$, and $H(a)$, $H_{-}(a)$ are Hilbert spaces of spinor fields defined in the subsections below. We will prove the\\

\noindent
\textbf{Proposition.} {\it For every $\sigma \in IKS(\Sigma)$ there exists a unique $\xi_{0}\in H(a)$ (resp. $\in H_{-}(a)$ if $M$ has a boundary) such that
$$\xi=f\mathcal{A}\sigma +\xi_{0} \in \text{Ker}\widehat{\mathfrak{D}} \ (\text{resp.} \ \in\text{Ker}\widehat{\mathfrak{D}} \cap H_{-}(a)) \quad \text{and} \quad \mathcal{H}(V_{\sigma},\alpha_{\sigma})\geq 0,$$
where $V_{\sigma }=\left\langle \sigma ,\sigma \right\rangle$ and $\alpha _{\sigma }(X)=\left\langle X \cdot e_{0}\cdot\sigma ,\sigma \right\rangle$ .}\\

\noindent
In section 3.3, it will proved that the couple $(V_{\sigma},\alpha_{\sigma})$ belongs to $N_{b}\oplus \mathfrak{Kill}(\Bbb H^{n})$ so that $\mathcal{H}(V_{\sigma},\alpha_{\sigma})$ is actually well-defined.\\
The computations we will make in section 3.1 prove that if we integrate the Bochner-Lichnerowicz-Weitzenböck-Witten formula with an asymptotically imaginary Killing spinor $f\mathcal{A}\sigma$, then the boundary integrals tend to some {\it global charge} $\mathcal{H}(V_{\sigma},\alpha_{\sigma})$, when $r$ goes to infinity. In fact it is still true if we {\it perturb}  $f\mathcal{A}\sigma$ with a smooth compactly supported spinor field $\xi _{0}$ (that is to say if we consider $f\mathcal{A}\sigma+\xi_{0}$ instead of $f\mathcal{A}\sigma$). Actually we will show in section 3.2 that we can find a {\it perturbation} $\xi _{0}$   in a relevant Hilbert space such that $\xi _{0}$ has no contribution at infinity, and $f \mathcal{A}\sigma + \xi_{0}$ belongs to the kernel of $\widehat{\mathfrak{D}}$.\\
This will naturally imply the non-negativity of $\mathcal{H}(V_{\sigma},\alpha_{\sigma})$ when $\sigma $ is a {\it $\beta $-imaginary Killing spinor}. This is the reason why we focus on the study of the Killing equation in section 3.3 so as to interpret the non-negativity of the $\mathcal{H}(V_{\sigma},\alpha_{\sigma})$.

\subsection{Energy-Momentum and Imaginary Killing Spinors}

The aim of this section is to show the\\

\noindent
\textbf{Proposition.} \textit{Let} $\xi =f \mathcal{A}\sigma + \xi_{0} $ \textit{, where} $\sigma\in IKS(\Sigma)$ \textit{and} $\xi_{0}$ \textit{is a compactly supported spinor field. Then we have}
\begin{eqnarray*}
\mathcal{H}(V_{\sigma},\alpha_{\sigma})&=& 4 \lim_{r\rightarrow+\infty} \int_{S_{r}}\left\langle \widehat{\nabla}_{\mathcal{A}\nu_{r}}\xi+ \mathcal{A}\nu_{r}\cdot\widehat{\mathfrak{D}}\xi, \xi \right\rangle_{\gamma}\\
&=&4\int_{M}\left(\left|\widehat{\nabla}\xi\right|^{2}_{\gamma}+\left\langle\widehat{\mathfrak{R}}\xi,\xi\right\rangle_{\gamma}\right)-4\int_{M} \left|\widehat{\mathfrak{D}}\xi\right|^{2}_{\gamma}.
\end{eqnarray*}

\noindent
Remark that the only important data is the exact {\it $\beta $-imaginary Killing spinor} $\sigma $ involved in the definition of the couple $(V_{\sigma},\alpha_{\sigma})$.\\

\noindent
{\sc Proof.}
Remember that
$$\int_{M_{r}} \left|\widehat{\mathfrak{D}}\psi\right|^{2}_{\gamma} = 
\int_{M_{r}}\left(\left|\widehat{\nabla}\psi\right|^{2}_{\gamma}+\left\langle \widehat{\mathfrak{R}}\psi,\psi\right\rangle_{\gamma}\right)
-\int_{S_{r}}\left\langle \widehat{\nabla}_{\mathcal{A}\nu_{r}}\psi+ \mathcal{A}\nu_{r}\cdot_{\gamma}\widehat{\mathfrak{D}}\psi, \psi\right\rangle_{\gamma} \text{dVol}_{S_{r}}  \quad ,
$$
where $\nu_{r}$ denotes the $b$-normal of $S_{r}$, $e_{0}=\partial_{t}$ and we set $e_{1}= \mathcal{A}\nu_{r}$ for the remainder of the proof.\\
We have to work on the expression $\left\langle \widehat{\nabla}_{\mathcal{A}\nu_{r}}\psi+ \mathcal{A}\nu_{r}\cdot_{\gamma}\widehat{\mathfrak{D}}\psi, \psi\right\rangle_{\gamma}$ in order to identify the integrand used to compute the energy-momentum for some couple $(f,\alpha )$. We start with noticing that  $e_{1}\cdot_{\gamma}e_{1}\cdot_{\gamma}\widehat{\nabla}_{e_{1}}= -\widehat{\nabla}_{e_{1}}=-\widehat{\nabla}_{\mathcal{A}\nu_{r}}$ so that
$$\qquad \widehat{\nabla}_{\mathcal{A}\nu_{r}}\psi+ \mathcal{A}\nu_{r}\cdot_{\gamma}\widehat{\mathfrak{D}}\psi= \mathcal{A}\nu_{r} \cdot_{\gamma} \left(\sum^{n}_{j=2}e_{j}\cdot_{\gamma}\widehat{\nabla}_{e_{j}}\right)\psi .$$
From now on we work on 
$$\left\langle \mathcal{A}\nu_{r} \cdot_{\gamma} \left(\sum^{n}_{j=2}e_{j}\cdot_{\gamma}\widehat{\nabla}_{e_{j}}\right)\cdot_{\gamma}\varphi,\varphi\right\rangle_{\gamma}.$$
Let us take $\sigma$ a {\it $\beta$-imaginary Killing spinor}, that is to say a spinor field solution, by definition, of $\widehat{D}_{X}\sigma = D_{X}\sigma + \frac{\textbf{i}}{2}X \cdot_{\beta} \sigma =0$,  for every vector field $X\in \Gamma (TM)$. Consider $f$ a smooth cutoff function which is 0 on $M$ except on a compact neighbourhood of the infinity boundary of $M$ where $f\equiv 1$. Then we have
$$\begin{array}{lllll}
\widehat{\nabla}_{X}(f \mathcal{A}\sigma)&=& \text{d}f(X)\mathcal{A}\sigma & +	f \widehat{\nabla}_{X}(\mathcal{A}\sigma)&\\
	&=& \text{d}f(X)\mathcal{A}\sigma & +f(\overline{\nabla}_{X}- \widetilde{\nabla}_{X})(\mathcal{A}\sigma) & + f\left(\widetilde{\nabla}_{X} +\frac{\textbf{i}}{2}  X\cdot_{\gamma} -\frac{1}{2}k(X) \cdot_{\gamma}e_{0}\cdot_{\gamma}\right)(\mathcal{A}\sigma),	
\end{array}$$
but since $\widetilde{\nabla}_{X}(\mathcal{A}\sigma)=\mathcal{A} \overline{D}_{X}\sigma=-\frac{\textbf{i}}{2} \mathcal{A}(X \cdot_{\beta}\sigma)=-\frac{i}{2} (\mathcal{A}X)\cdot_{\gamma} (\mathcal{A}\sigma)$, we obtain
$$\widehat{\nabla}_{X}(f\mathcal{A}\sigma)= \text{d}f(X)\mathcal{A}\sigma  +f(\overline{\nabla}_{X}- \widetilde{\nabla}_{X})(\mathcal{A}\sigma) -\frac{1}{2} f\left( k(X)\cdot_{\gamma}e_{0}+\textbf{i}(\mathcal{A}-Id)X \right) \cdot_{\gamma}(\mathcal{A}\sigma),$$
that we restrict to the neighbourhood where $f\equiv 1$
$$\widehat{\nabla}_{X}(\mathcal{A}\sigma)= (\overline{\nabla}_{X}- \widetilde{\nabla}_{X})(\mathcal{A}\sigma) -\frac{1}{2} (k(X)\cdot_{\gamma}e_{0}+\textbf{i}(\mathcal{A}-Id)X) \cdot_{\gamma}(\mathcal{A}\sigma).$$
As a consequence our boundary term becomes for $r$ great enough 
$$\sum^{n}_{j=2}\left\langle  \mathcal{A}\nu_{r} \cdot_{\gamma} e_{j}\cdot_{\gamma}\left( (\overline{\nabla}_{e_{j}}- \widetilde{\nabla}_{e_{j}})-\frac{1}{2} (k(e_{j})\cdot_{\gamma}e_{0}+\textbf{i}(\mathcal{A}-Id)e_{j}) \cdot_{\gamma} \right) (\mathcal{A}\sigma),\mathcal{A}\sigma \right\rangle_{\gamma}.$$
We will estimate this boundary term in several steps. From the decay assumptions stated section in 1.2, the gauge is supposed to be of the form $\mathcal{A}=Id + B + O(|B|^{2})$, where \textit{B} has the same decay to 0 as $e=g-b$. In the following $(\epsilon_{k}=\mathcal{A}^{-1}e_{j})^{n}_{j=0}$ is a $\beta$-orthonormal frame.\\

\noindent
We begin with the easiest term
\begin{eqnarray*}
	\sum^{n}_{j=2}\left\langle  \mathcal{A}\nu_{r} \cdot_{\gamma} e_{j}\cdot_{\gamma}k(e_{j})\cdot_{\gamma}e_{0}\cdot_{\gamma}(\mathcal{A}\sigma), \mathcal{A}\sigma \right\rangle_{\gamma}&=&\sum^{n}_{j=2} \left\langle  \mathcal{A}\nu_{r} \cdot_{\gamma} \mathcal{A}\epsilon_{j}\cdot_{\gamma}k(\mathcal{A}\epsilon_{j})\cdot_{\gamma}\mathcal{A}\epsilon_{0}\cdot_{\gamma}(\mathcal{A}\sigma), \mathcal{A}\sigma \right\rangle_{\gamma}\\
	&=& \sum^{n}_{j=2}\left\langle  \nu_{r} \cdot_{\beta} \epsilon_{j}\cdot_{\beta}\mathcal{A}^{-1}\circ k\circ \mathcal{A}(\epsilon_{j})\cdot_{\beta}\epsilon_{0}\cdot_{\beta}\sigma, \sigma \right\rangle_{\beta}.	
\end{eqnarray*}
But we note that  $\mathcal{A}^{-1}\circ k\circ \mathcal{A}= k -B \circ k+k\circ B + O(|B|^{2})$. Now \textit{B} has the same decay as $k$ so $B \circ k+k\circ B = O(|B|^{2})$, terms that we can neglect since the energy-momentum is computed by a limit procedure of integrals over large spheres. We conclude that $\mathcal{A}^{-1}\circ k \circ \mathcal{A} \approx k$, where for convenience the relation  $\diamond  \approx  \star $ means that $|\diamond -\star   |$ is at least a $O(e^{-2\tau r})$ when $r$ goes to infinity. Moreover 
\begin{eqnarray*}
	\nu_{r} \cdot_{\beta} \sum^{n}_{j=2} \epsilon_{j}\cdot_{\beta} k(\epsilon_{j}) &=&\epsilon_{1} \cdot_{\beta} \left(\sum^{n}_{j=1} \epsilon_{j}\cdot_{\beta} k(\epsilon_{j})- \epsilon_{1}\cdot_{\beta} k(\epsilon_{1}) \right)\\
	&=& k(\nu_{r})-(\text{tr}_{b}k)\nu_{r},
\end{eqnarray*}
which implies
\begin{eqnarray*}
	\sum^{n}_{j=2}\left\langle  \mathcal{A}\nu_{r} \cdot_{\gamma} e_{j}\cdot_{\gamma}k(e_{j})\cdot_{\gamma}e_{0}\cdot_{\gamma}(\mathcal{A}\sigma), \mathcal{A}\sigma \right\rangle_{\gamma}&\approx & \left\langle  k(\nu_{r})-(\text{tr}_{b}k)\nu_{r} \cdot_{\beta} \epsilon_{0}\cdot_{\beta} \sigma,\sigma \right\rangle_{\beta}	\\
	&= & \left(i_{\alpha_{\sigma}}k-(\text{tr}_{b}k)\alpha_{\sigma}\right)(\nu_{r}),
\end{eqnarray*}
where $\alpha_{\sigma}(X)= \left\langle X \cdot_{\beta} \epsilon_{0} \cdot_{\beta} \sigma,\sigma\right\rangle_{\beta}$. \\

\noindent
The second term we study is
\begin{eqnarray*}
\textbf{i}\sum^{n}_{j=2}	\left\langle  \mathcal{A}\nu_{r} \cdot_{\gamma} e_{j}\cdot_{\gamma}(\mathcal{A}-Id)(e_{j})\cdot_{\gamma}(\mathcal{A}\sigma), \mathcal{A}\sigma \right\rangle_{\gamma}&=& \textbf{i}\sum^{n}_{j=2}\left\langle  \mathcal{A}\nu_{r} \cdot_{\gamma} \mathcal{A}\epsilon_{j}\cdot_{\gamma}(\mathcal{A}-Id)(\mathcal{A}\epsilon_{j})\cdot_{\gamma}(\mathcal{A}\sigma), \mathcal{A}\sigma \right\rangle_{\gamma}\\
	&=&\textbf{i}\sum^{n}_{j=2} \left\langle  \nu_{r} \cdot_{\beta} \epsilon_{j}\cdot_{\beta}\mathcal{A}^{-1}\circ (\mathcal{A}-Id)\circ \mathcal{A}(\epsilon_{j})\cdot_{\beta}\sigma, \sigma \right\rangle_{\beta}\\
	&\approx & \textbf{i}\sum^{n}_{j=2}\left\langle  \nu_{r} \cdot_{\beta} \epsilon_{j}\cdot_{\beta}B(\epsilon_{j})\cdot_{\beta}\sigma, \sigma \right\rangle_{\beta},
\end{eqnarray*}
but thanks to the same property as above 
$$\nu_{r}\cdot_{\beta}\sum^{n}_{j=2}\epsilon_{j}\cdot_{\beta}B(\epsilon_{j})=B(\nu_{r})-(\text{tr}_{b}B)\nu_{r},$$
which induces
\begin{eqnarray*}
	\textbf{i}\sum^{n}_{j=2}\left\langle  \mathcal{A}\nu_{r} \cdot_{\gamma} e_{j}\cdot_{\gamma}(\mathcal{A}-Id)(e_{j})\cdot_{\gamma}(\mathcal{A}\sigma), \mathcal{A}\sigma \right\rangle_{\gamma}&\approx &\textbf{i} \left\langle  B(\nu_{r})-(\text{tr}_{b}B)\nu_{r} \cdot_{\beta} \sigma,\sigma \right\rangle_{\beta}\\
&=&\left(i_{\nabla V_{\sigma}}B-(\text{tr}_{b}B)\text{d}V_{\sigma}\right)(\nu_{r}),
\end{eqnarray*}
where $\text{d}V_{\sigma}(X)= \textbf{i}\left\langle X \cdot_{\beta} \sigma,\sigma\right\rangle_{\beta}$.\\

\noindent
The last term we have to study is certainly the most difficult (summation convention $k\in\{2,3, \cdots,n\},l\in\{1,2, \cdots,n\}, m\in\{1,2,\cdots,n\}$)
$$\begin{array}{lll}
	\left\langle  \mathcal{A}\nu_{r} \cdot_{\gamma} e_{k}\cdot_{\gamma} (\overline{\nabla}_{e_{k}}- \widetilde{\nabla}_{e_{k}}) (\mathcal{A}\sigma),\mathcal{A}\sigma \right\rangle_{\gamma} &=& 
\frac{1}{4}\left\langle (\overline{\omega}_{lm}-\widetilde{\omega}_{lm})(e_{k})  \mathcal{A}\nu_{r} \cdot_{\gamma} e_{k}\cdot_{\gamma} e_{l}\cdot_{\gamma}e_{m}\cdot_{\gamma}(\mathcal{A}\sigma),\mathcal{A}\sigma\right\rangle_{\gamma}\\
&=&\frac{1}{4} \left\langle (\overline{\omega}_{lm}-\widetilde{\omega}_{lm})\circ \mathcal{A} (\epsilon_{k}) \nu_{r} \cdot_{\beta} \epsilon_{k} \cdot_{\beta}\epsilon_{l}\cdot_{\beta} \epsilon_{m}\cdot_{\beta}\sigma , \sigma \right\rangle_{\beta}\\
&=& \frac{1}{4} S
\end{array}$$

\begin{eqnarray*}
	S &=&   \sum^{n}_{k,l,m=2}\left\langle (\overline{\omega}_{lm}-\widetilde{\omega}_{lm})(e_{k})  \mathcal{A}\nu_{r} \cdot_{\gamma} e_{k}\cdot_{\gamma} e_{l}\cdot_{\gamma}e_{m}\cdot_{\gamma}(\mathcal{A}\sigma),\mathcal{A}\sigma\right\rangle_{\gamma}\\
	& +& 2\sum^{n}_{k,l=2} \left\langle (\overline{\omega}_{1l}-\widetilde{\omega}_{1l})(e_{k})  \mathcal{A}\nu_{r} \cdot_{\gamma} e_{k}\cdot_{\gamma} e_{1}\cdot_{\gamma}e_{l}\cdot_{\gamma}(\mathcal{A}\sigma),\mathcal{A}\sigma\right\rangle_{\gamma}\\
	&=& S_{1}+2S_{2}	.
\end{eqnarray*}
We will give estimates of each $S_{k}$, keeping in mind that they are real and that every term that is  at least $O(|B|^{2})$ can be neglected when $r \rightarrow  +\infty$ for the computations of the {\it global charge integrals}.\\

\textbf{Estimate of $S_{1}$}

$$S_{1}=\sum^{n}_{k,l,m=2}\left\langle (\overline{\omega}_{lm}-\widetilde{\omega}_{lm})(e_{k})  \mathcal{A} \nu_{r} \cdot_{\gamma} e_{k}\cdot_{\gamma} e_{l}\cdot_{\gamma}e_{m}\cdot_{\gamma}(\mathcal{A}\sigma),\mathcal{A}\sigma\right\rangle_{\gamma}.$$
We can keep only the subscripts $l\neq m$ because of the skew-symmetry of $(\omega-\widetilde{\omega})$. Besides if we suppose that $k=l$, we have terms like $\left\langle \mathcal{A}\nu_{r} \cdot_{\gamma} e_{m}\cdot_{\gamma}(\mathcal{A}\sigma),\mathcal{A}\sigma\right\rangle_{\gamma}$ which belong to $\textbf{i}\mathbb{R}$. So we can sum over $k,l, m$ distinct subscripts without any loss of generality. On the other hand
$$ (\overline{\omega}_{lm}-\widetilde{\omega}_{lm})(e_{k})= \frac{1}{2} \left(-g(\widetilde{T}(e_{k},e_{l}),e_{m})+g(\widetilde{T}(e_{k},e_{m}),e_{l})+g(\widetilde{T}(e_{l},e_{m}),e_{k}) \right)$$
where the two last terms of the right-hand side member are symmetric with respect to $(l,k)$, so they vanish when we sum over $k$ and $l$ distinct. Consequently
$$\begin{array}{lll}
(\overline{\omega}_{lm}-\widetilde{\omega}_{lm})(e_{k})   \epsilon_{k}\cdot_{\beta} \epsilon_{l}\cdot_{\beta}\epsilon_{m} &=&
\frac{1}{2} b\left( \mathcal{A}^{-1}(\overline{D}_{e_{k}}\mathcal{A})\epsilon_{l}-\mathcal{A}^{-1}(\overline{D}_{e_{l}}\mathcal{A})\epsilon_{k},\epsilon_{m}\right)\epsilon_{k}\cdot_{\beta} \epsilon_{l}\cdot_{\beta}\epsilon_{m}\\
&=&  b\left( \mathcal{A}^{-1}(\overline{D}_{e_{k}}\mathcal{A})\epsilon_{l},\epsilon_{m}\right) \epsilon_{k}\cdot_{\beta} \epsilon_{l}\cdot_{\beta}\epsilon_{m},
\end{array}$$
but
$$\begin{array}{lll}
b( \mathcal{A}^{-1}(\overline{D}_{e_{k}}\mathcal{A})\epsilon_{l},\epsilon_{m})&=& b( \mathcal{A}^{-1}(\overline{D}_{e_{k}}(\mathcal{A}\epsilon_{l})-\mathcal{A}\overline{D}_{e_{k}}\epsilon_{l}),\epsilon_{m})\\
&\approx& b( \overline{D}_{\epsilon_{k}}(B\epsilon_{l})- B(\overline{D}_{\epsilon_{k}}\epsilon_{l}),\epsilon_{m}) \\
&=& b( (\overline{D}_{\epsilon_{k}}B)\epsilon_{l},\epsilon_{m}) ,	
\end{array}$$
expression which is symmetric with respect to $(l,m)$, since $\overline{D}B$ is a symmetric endomorphism. Consequently
$$\sum_{\text{\textit{k,l,m} distinct}}
(\overline{\omega}_{lm}-\widetilde{\omega}_{lm})(e_{k})   \epsilon_{k}\cdot_{\beta} \epsilon_{l}\cdot_{\beta}\epsilon_{m} \approx 0 ,$$
when  $r\rightarrow+\infty$.\\

\textbf{Estimate of $S_{2}$}

\begin{eqnarray*}
	S_{2}&=& \sum^{n}_{k,l=2} \left\langle (\overline{\omega}_{1l}-\widetilde{\omega}_{1l})(e_{k})  \mathcal{A}\nu_{r} \cdot_{\gamma} e_{k}\cdot_{\gamma} e_{1}\cdot_{\gamma}e_{l}\cdot_{\gamma}(\mathcal{A}\sigma),\mathcal{A}\sigma\right\rangle_{\gamma}\\
	&=& -  \sum^{n}_{k=2}\left\langle (\overline{\omega}_{1k}-\widetilde{\omega}_{1k})(e_{k})\sigma,\sigma\right\rangle_{\beta}\\
	&& +\sum_{k\neq l}\left\langle (\overline{\omega}_{1l}-\widetilde{\omega}_{1l})(e_{k}) e_{k}\cdot_{\gamma}e_{l}\cdot_{\gamma}(\mathcal{A}\sigma),\mathcal{A}\sigma\right\rangle_{\gamma},
\end{eqnarray*}
but the second sum is in $\textbf{i}\mathbb{R}$, so it remains
$$ \Re e (S_{2})= - \left( \sum^{n}_{k=1}(\overline{\omega}_{1k}-\widetilde{\omega}_{1k})(e_{k}) \right)V_{\sigma}.$$
We only have to compute
\begin{eqnarray*}
-\sum^{n}_{k=1}(\overline{\omega}_{1k}-\widetilde{\omega}_{1k})(e_{k})&=& \sum^{n}_{k=1} g( (\overline{D}_{e_{1}}\mathcal{A})\mathcal{A}^{-1}e_{k},e_{k}-(\overline{D}_{e_{k}}\mathcal{A})\mathcal{A}^{-1}e_{1},e_{k})\\
&=& S'_{2}-S''_{2}
\end{eqnarray*}
We focus on
\begin{eqnarray*}
S''_{2}  &=& \sum^{n}_{k=1}g( (\overline{D}_{e_{k}}\mathcal{A})\mathcal{A}^{-1}e_{1},e_{k})\\
        &=&\sum^{n}_{k=1}b( \mathcal{A}^{-1}(\overline{D}_{e_{k}}\mathcal{A})\epsilon_{1},\epsilon_{k})\\
	&=& \sum^{n}_{k=1}b( \mathcal{A}^{-
	1}\overline{D}_{e_{k}}(\mathcal{A}\epsilon_{1})-\overline{D}_{e_{k}}\epsilon_{1},\epsilon_{k})\\
	&\approx &\sum^{n}_{k=1}b( (\overline{D}_{e_{k}}B)\epsilon_{1},\epsilon_{k})\\
	&\approx&\sum^{n}_{k=1}b(\epsilon_{1},(\overline{D}_{\epsilon_{k}}B)\epsilon_{k})\\
	&=&-\text{div}_{b}B(\nu_{r}).
\end{eqnarray*}
As regards the first term $ S_{2}'$, we decompose the gauge endomorphism $\mathcal{A}$ as follows:${\mathcal{A}\epsilon_{i}=\sum^{n}_{k=0}\mathcal{A}^{k}_{i}\epsilon_{k}}$. We remind that  $\mathcal{A}\epsilon_{0}=\epsilon_{0}$, $\mathcal{A}(TM)\subset TM$ and so we have $\mathcal{A}^{k}_{0}=\mathcal{A}^{0}_{k}=0, \  k \geq1$.

\begin{eqnarray*}
 S_{2}'&=&\sum^{n}_{k=1}g( (\overline{D}_{e_{1}}\mathcal{A})\mathcal{A}^{-1}e_{k},e_{k})\\
 &=& \sum^{n}_{k=1}b( (\mathcal{A}^{-1}\overline{D}_{e_{1}}\mathcal{A})\epsilon_{k},\epsilon_{k})\\
 &=& \sum^{n}_{k=1} b(\mathcal{A}^{-1}\overline{D}_{e_{1}}(\mathcal{A}\epsilon_{k})-\overline{D}_{e_{1}}\epsilon_{k},\epsilon_{k})\\
 &\approx& \sum^{n}_{k,l=1} (e_{1}\cdot \mathcal{A}^{l}_{k})  \left\{b( \epsilon_{l},\epsilon_{k}) -b( B\epsilon_{l},\epsilon_{k})\right\} -\sum^{n}_{k=1}b( \overline{D}_{e_{1}}\epsilon_{k},\epsilon_{k})\\
 &&+\sum^{n}_{k,l=1}  \mathcal{A}^{l}_{k}\left\{b( \overline{D}_{e_{1}}\epsilon_{l},\epsilon_{k}) -b( B\overline{D}_{e_{1}}\epsilon_{l},\epsilon_{k})\right\}\\
 &\approx& (\epsilon_{1}\cdot\text{tr}_{b}B)-\sum^{n}_{k=1}b( \overline{D}_{e_{1}}\epsilon_{k},\epsilon_{k}) +\sum^{n}_{k,l=1} (\delta^{k}_{l}+ B^{l}_{k})\left\{b( \overline{D}_{e_{1}}\epsilon_{l},\epsilon_{k}) -b( B\overline{D}_{e_{1}}\epsilon_{l},\epsilon_{k})\right\}\\
 &\approx&(\epsilon_{1}\cdot\text{tr}_{b}B) -\sum^{n}_{k=1}\left\{b( B\overline{D}_{e_{1}}\epsilon_{k},\epsilon_{k}) -b( \overline{D}_{e_{1}}\epsilon_{k},B\epsilon_{k})\right\}\\
 &=& \text{d}(\text{tr}_{b}B)(\nu_{r}),
\end{eqnarray*} 
that entails
$$\Re e (S_{2})\approx V_{\sigma}(\text{d}(\text{tr}_{b}B)+\text{div}_{b}B)(\nu_{r}).$$\\

\noindent
We can conclude, taking $B=-\frac{1}{2}e$, that the real part of our boundary integrand is nothing but
$$\frac{1}{4}\left(-V_{\sigma}(\delta_{b}e + \text{d}tr_{b}e) - i_{\nabla^{b} V_{\sigma}}e +(tr_{b}e) \text{d}V_{\sigma} - 2i_{\alpha^{\sharp}_{\sigma}}k + 2(tr_{b}k)\alpha_{\sigma}\right) (\nu_{r}),$$
what achieves the proof. \hspace{\stretch{4}}$\square$

\subsection{Analysis of $\widehat{\mathfrak{D}}$}

\noindent
This section is devoted to the study of the analytical properties of $\widehat{\mathfrak{D}}$. The first paragraph deals with the case where $M$ has no boundary whereas the second one deals with the case where $M$ has a compact and connected boundary denoted as usual by $\partial M$.

\subsubsection{$M$ without Boundary}

\noindent
\textbf{Proposition.} {\it For every $\sigma {\in} IKS(\Sigma)$ there exists a unique $\xi_{0}{\in} H(a)$ such that $\xi{=}f\mathcal{A}\sigma {+}\xi_{0} {\in} \text{Ker} \widehat{\mathfrak{D}}$ and 
$$\mathcal{H}(V_{\sigma},\alpha_{\sigma})= \lim_{r\rightarrow+\infty} \int_{S_{r}}\left\langle \widehat{\nabla}_{\mathcal{A}\nu_{r}}\xi+ \mathcal{A}\nu_{r}\cdot\widehat{\mathfrak{D}}\xi, \xi\right\rangle_{\gamma}= 4\int_{M}\left(\left|\widehat{\nabla}\xi\right|^{2}_{\gamma}+\left\langle\widehat{\mathfrak{R}}\xi,\xi\right\rangle_{\gamma}\right)\geq 0.$$ }\\

\noindent
{\sc Proof.}
We study in a usual way the analytical properties of $\widehat{\mathfrak{D}}$. Let us consider ${C^{\infty}_{0}(\Sigma)=C^{\infty}_{0}}$ the space of smooth and compactly supported spinors. We define a sesquilinear form on $C^{\infty}_{0}$ by
$$ a(\varphi,\psi)= \int_{M}\left\langle \widehat{\mathfrak{D}}\varphi, \widehat{\mathfrak{D}}\psi\right\rangle_{\gamma}\text{d}\mu_{g},$$
where $\text{d}\mu_{g}$ denotes the standard volum form of the metric g. The form $a$ is clearly bounded and non-negative on $C^{\infty}_{0}$. We define the usual Sobolev space
$$H^{1}(\Sigma)=\left\{\psi\in\Sigma / \int_{M}\left|\psi\right|^{2}_{\gamma}+ \left|\nabla \psi\right|^{2}_{\gamma}<\infty \right\}.$$\\

\noindent
\textbf{Definition.} \textit{We set} $H(a):=  \overline{C^{\infty}_{0}}^{a}.$\\

\noindent
\textbf{Remark. Weighted Poincar\'e inequality}\\
$$ \exists\omega \in L_{loc}^{1}(M,\text{dVol}_{g})\quad ess_{M}inf \omega > 0\quad \forall u\in C_{0}^{1} \quad \int_{M}\omega|u|^{2}\text{dVol}_{g}\leq \int_{M}|\widehat{\nabla}u|^{2}\text{dVol}_{g}.$$
It is easy to see that $\Gamma$, the symmetric part of the connection $\widehat{\nabla}$ is given by $\Gamma_{X}=\frac{1}{2}\left\{ k(X)\cdot e_{0}-\textbf{i}X\right\}\cdot$, and so satisfies the conditions (cf. \cite{ChB}) in order to have the existence of a  weighted Poincar\'e inequality that is to say $\Gamma\in L^{n}_{loc}(M) \quad \text{and} \quad \limsup_{x \rightarrow 0 }|x\Gamma_{x}|< \frac{n-1}{2}$. Such a weighted Poincaré inequality insures the continuity of the embedding of $\xymatrix@1{H(a)\ar@{^{(}->}[r]& H^{1}_{loc}}$ (cf. section 3.2.2 for a proof of this fact when $M$ has a compact boundary).\hspace{\stretch{4}}$\square$\\

%\noindent
%Since for every $\psi \in C^{\infty}_{0}$ we have
%$$\left(1+\frac{n}{4} \right)^{-1}||\psi||^{2}_{1} \leq \int_{M} \left|\widehat{\mathfrak{D}}\psi\right|^{2}_{\gamma} = 
%\int_{M}\left(\left|\nabla\psi\right|^{2}_{\gamma}+\frac{n}{4}\left|\psi\right|^{2}+\left\langle \widehat{\mathfrak{R}}\psi,\psi\right\rangle_{\gamma}\right)
% \quad ,
%$$
%we can claim that $\xymatrix@1{H(a)\ar@{^{(}->}[r]& H^{1}}$ is a continuous embedding. 

\noindent
We notice that for $r$ great enough and for $\sigma \in IKS(\Sigma)$
\begin{eqnarray*}
	\widehat{\nabla}_{X}(\mathcal{A}\sigma)&=& \nabla_{X}(\mathcal{A}\sigma) + \frac{\textbf{i}}{2}X \cdot_{\gamma}(\mathcal{A}\sigma)\\
	&=& (\overline{\nabla}_{X}- \widetilde{\nabla}_{X})(\mathcal{A}\sigma) -\frac{1}{2} (k(X)\cdot_{\gamma}e_{0}+\textbf{i}(\mathcal{A}-Id)X) \cdot_{\gamma}(\mathcal{A}\sigma).\\
\end{eqnarray*}
But the relations
$$\left\{\begin{array}{rll}
	\widetilde{T}(X,Y)&=&-((\overline{D}_{X}\mathcal{A})\mathcal{A}^{-1}Y-(\overline{D}_{Y}\mathcal{A})\mathcal{A}^{-1}X)\\
	2g( \widetilde{\nabla}_{X}Y-\overline{\nabla}_{X}Y,Z )&=&g(\widetilde{T}(X,Y),Z )-g(\widetilde{T}(X,Z),Y )-g(\widetilde{T}(Y,Z),X)
\end{array} \right.$$
tell us that $|(\overline{\omega}_{ij} - \widetilde{\omega_{ij}})(e_{k})|\leq C |\mathcal{A}^{-1}||\overline{D}\mathcal{A}|$. We get an estimate
$$ |\widehat{\mathfrak{D}}(\mathcal{A}\sigma)|\leq C|\mathcal{A}|(|\overline{D}\mathcal{A}|+|\mathcal{A}-Id|+|k|)|\sigma|\in L^{2}(M,\text{d}\mu_{g}),$$
which infers that $\widehat{\mathfrak{D}}(f\mathcal{A}\sigma)\in L^{2}(M,\text{d}\mu_{g})$. We now consider the linear form \textit{l} on $H(a)$ defined by
$$l(\psi)= \int_{M}\left\langle\widehat{\mathfrak{D}}(f\mathcal{A}\sigma),\widehat{\mathfrak{D}}\psi\right\rangle_{\gamma}\text{d}\mu_{g}.$$ 
Thanks to our estimate above we get $|l(\psi)|^{2} \leq \left\|\widehat{\mathfrak{D}}(f\mathcal{A}\sigma)\right\|^{2}_{L^{2}}a(\psi,\psi)$, that gives the continuity of \textit{l} in $H(a)$. We can claim, thanks to Lax-Milgram theorem, that there exists a unique $\xi_{0}\in H(a)$ such that $l=a(-\xi_{0},\cdot)$. In other words
$$\int_{M}\left\langle(\widehat{\mathfrak{D}})^{*}\widehat{\mathfrak{D}}(f\mathcal{A}\sigma+\xi_{0}),\psi\right\rangle_{\gamma}=0.$$
Since $\widehat{\mathfrak{D}}^{*}=\widehat{\mathfrak{D}} +\textbf{i}n$, we have in the distributional sense $(\widehat{\mathfrak{D}} +\textbf{i}n)\widehat{\mathfrak{D}}\xi =0$, where we have set $\xi= f\mathcal{A}\sigma+\xi_{0}$. By an elliptic regularity argument, $\widehat{\mathfrak{D}}\xi$ is in fact smooth and $(\widehat{\mathfrak{D}})^{k}\xi$ are $L^{2}$, for every $k\in\mathbb{N}$. It follows
\begin{eqnarray*}
\int_{M}\left\langle(\widehat{\mathfrak{D}})^{2}\xi,(\widehat{\mathfrak{D}})^{2}\xi\right\rangle_{\gamma} &=& \int_{M}\left\langle(\widehat{\mathfrak{D}} +\textbf{i}n)(\widehat{\mathfrak{D}})^{2}\xi,\widehat{\mathfrak{D}}\xi\right\rangle_{\gamma}\\
&=&\int_{M}\left\langle \widehat{\mathfrak{D}}(\widehat{\mathfrak{D}} +\textbf{i}n)\widehat{\mathfrak{D}}\xi,\widehat{\mathfrak{D}}\xi\right\rangle_{\gamma}\\
&=& 0,
\end{eqnarray*}
that implies $(\widehat{\mathfrak{D}})^{2}\xi=0$, but we already know that $(\widehat{\mathfrak{D}} +\textbf{i}n)\widehat{\mathfrak{D}}\xi =0$, and thereby $ \widehat{\mathfrak{D}}\xi=0$. We now apply our integration formula to $\xi$
\begin{eqnarray*}
\mathcal{H}(V_{\sigma},\alpha_{\sigma})&=&  \lim_{r\rightarrow+\infty} \int_{S_{r}}\left\langle \widehat{\nabla}_{\mathcal{A}\nu_{r}}\xi+ \mathcal{A}\nu_{r}\cdot\widehat{\mathfrak{D}}\xi, \xi\right\rangle_{\gamma}\\
&=&4\int_{M}\left(\left|\widehat{\nabla}\xi\right|^{2}_{\gamma}+\left\langle\widehat{\mathfrak{R}}\xi,\xi\right\rangle_{\gamma}\right)-4\int_{M} \left|\widehat{\mathfrak{D}}\xi\right|^{2}_{\gamma}\\
&=&4\int_{M}\left(\left|\widehat{\nabla}\xi\right|^{2}_{\gamma}+\left\langle\widehat{\mathfrak{R}}\xi,\xi\right\rangle_{\gamma}\right)\geq 0,
\end{eqnarray*}
and the proof is complete.\hspace{\stretch{4}}$\square$

\subsubsection{$M$ with Boundary}

We will consider, in this section, a Riemannian slice $M$ that has a non empty inner boundary $\partial M$. $\breve{g},\stackrel{\smile}{\nabla},\breve{k}$ will denote respectively the induced metric, the connection and the second fundamental form which is defined by
\begin{eqnarray*}
  \overline{\nabla}_{X}Y&=&\stackrel{\smile}{\nabla}_{X}Y -\breve{k}(X,Y)\nu\\
  \overline{\nabla}_{X}\psi&=&\stackrel{\smile}{\nabla}_{X}\psi-\frac{1}{2}\breve{k}(X)\cdot\nu\cdot\psi,\\
\end{eqnarray*}
where $\nu$ is the normal to $\partial M$ pointing toward infinity (that is to say pointing inside), and $\cdot$ still denotes the Clifford action with respect to the metric $\gamma$. Consequently our integration formula has another boundary term
$$\int_{M_{r}} \left|\widehat{\mathfrak{D}}\psi\right|^{2}_{\gamma} = 
\int_{M_{r}}\left(\left|\widehat{\nabla}\psi\right|^{2}_{\gamma}+\left\langle \widehat{\mathfrak{R}}\psi,\psi\right\rangle_{\gamma}\right)
-\int_{S_{r}}\left\langle \widehat{\nabla}_{\mathcal{A}\nu_{r}}\psi+ \mathcal{A}\nu_{r}\cdot\widehat{\mathfrak{D}}\psi, \psi\right\rangle_{\gamma} + \int_{\partial M}\left\langle \widehat{\nabla}_{\nu}\psi+ \nu\cdot\widehat{\mathfrak{D}}\psi, \psi\right\rangle_{\gamma} .
$$
But if $\psi$ is a compactly supported smooth spinor field then, making $r \rightarrow \infty$ one finds
$$\int_{M} \left|\widehat{\mathfrak{D}}\psi\right|^{2}_{\gamma} = 
\int_{M}\left(\left|\widehat{\nabla}\psi\right|^{2}_{\gamma}+\left\langle \widehat{\mathfrak{R}}\psi,\psi\right\rangle_{\gamma}\right)
+ \int_{\partial M}\left\langle \widehat{\nabla}_{\nu}\psi+ \nu\cdot\widehat{\mathfrak{D}}\psi, \psi\right\rangle_{\gamma} .
$$
We then have to estimate the boundary integrand  $\left\langle \widehat{\nabla}_{\nu}\psi+ \nu\cdot\widehat{\mathfrak{D}}\psi,\psi \right\rangle$.\\

\noindent
\textbf{Lemma.} \textit{If} $(\nu=e_{1},e_{2},\cdots,e_{n})$ \textit{is a local orthonormal frame of} $TM_{|\partial M}$ \textit{then}
$$ \widehat{\nabla}_{\nu}\psi+ \nu\cdot\widehat{\mathfrak{D}}\psi = \nu \cdot \sum^{n}_{k=2}\widehat{\nabla}_{e_{k}}\psi \quad.$$\\

\noindent
\textsc{Proof.} Just remark that $\widehat{\nabla}_{\nu}\psi=-e_{1}\cdot e_{1}\cdot\widehat{\nabla}_{e_{1}}\psi $ \hspace{\stretch{4}}$\square$\\

\noindent
\textbf{Lemma.} \textit{Keeping our orthonormal frame} $(\nu=e_{1},e_{2},\cdots,e_{n}),$ \textit{we have}
$$\widehat{\nabla}_{\nu}\psi+ \nu\cdot\widehat{\mathfrak{D}}\psi =  \sum^{n}_{k=2}\nu \cdot e_{k}\cdot\stackrel{\smile}{\nabla}_{e_{k}}\psi+ \frac{1}{2}\left\{-\text{tr}\breve{k}-(n-1)\textbf{i}\nu+(\text{tr}k)\nu\cdot e_{0}-k(\nu)\cdot e_{0}\right\}\cdot \psi \quad.$$\\

\noindent
\textsc{Proof.} Using the formula above, we then express $\widehat{\nabla}$ in term of the $(n-1)$-dimensional connection and second form, and the $n$-dimensional second form.\hspace{\stretch{4}}$\square$\\

\noindent
Let us define $F \in End(\Sigma_{|\partial M})$  by $F(\psi)=\textbf{i}\nu\cdot \psi$. We sum up some basic properties of $F$ in the following\\

\noindent
\textbf{Proposition.}\textit{ The endomorphism F is symmetric, isometric with respect to} $\left\langle\cdot,\cdot\right\rangle$, \textit{commutes to the action of} $\nu\cdot$ \textit{ and anticommutes to each} $e_{k}\cdot,\, (k \neq 1).$\\

\noindent
\textbf{Lemma.} \textit{If} $F(\psi)=-\psi$ \textit{then}
$$\left\langle \widehat{\nabla}_{\nu}\psi+ \nu\cdot\widehat{\mathfrak{D}}\psi, \psi\right\rangle_{|\partial M}=\frac{1}{2}\left\langle e_{0}\cdot\left((-\text{tr}\breve{k}+(n-1))e_{0}+k(\nu)\right)\cdot \psi,\psi\right\rangle.$$\\

\noindent
\textsc{Proof.} Using the proposition above we know that $\nu\cdot e_{k} \, (k \neq 1)$ anticommutes with $F$ and the formula follows since $F$ respects  $\left\langle\cdot,\cdot\right\rangle$.\hspace{\stretch{4}}$\square$\\

\noindent
\textbf{Assumption.}  \textit{Let us suppose that the 4-vector} $\vec{k}:=(-\text{tr}\breve{k}+(n-1))e_{0}+k(\nu)$ \textit{is causal and positively oriented, that is to say} $\gamma(\vec{k},\vec{k})\leq0$ \textit{and} $\text{tr}\breve{k}\leq(n-1).$\\

\noindent
This assumption (which is exactly the same as for $\widehat{\mathfrak{R}}$) guarantees the non-negativity of the boundary integrand term $\left\langle \widehat{\nabla}_{\nu}\psi+ \nu\cdot\widehat{\mathfrak{D}}\psi, \psi\right\rangle_{|\partial M}=\frac{1}{2}\left\langle e_{0}\cdot \vec{k} \cdot \psi,\psi\right\rangle,$ whenever the boundary condition $F(\psi)=-\psi$ is satisfied. Although this assumption is vectorial, it clearly extends the one given in \cite{ChH}.\\
Let us define  $H_{-}(a)=\left\{\psi\in H(a)/ F(\psi)=-\psi\right\}$ where $H(a)$ has been defined in section 3.6.1. Still taking $\psi$ a compactly supported smooth spinor field in $H_{-}(a)$, we have 
$$a(\psi,\psi)=\int_{M}\left(\left|\widehat{\nabla}\psi\right|^{2}_{\gamma}+\left\langle \widehat{\mathfrak{R}}\psi,\psi\right\rangle_{\gamma}\right)+\frac{1}{2}\int_{\partial M}\left\langle e_{0}\cdot\vec{k}\cdot\psi,\psi\right\rangle,$$
whose each single term is non-negative tanks to our assumption.\\

\noindent
\textbf{Lemma.} $H_{-}(a)$ \textit{continuously embeds in} $H^{1}_{loc}$ \textit{and furthermore}
$$\left( \xymatrix@1{\psi_{k} \ar[r]^{H_{-}(a)} &\psi}\right) \Rightarrow  \left(
\xymatrix@1{\widehat{\nabla}\psi_{k} \ar[r]^{L^{2}(M)} & \widehat{\nabla}\psi} \quad \text{and}\quad \forall\Omega\subset M \quad |\Omega|<\infty \quad \xymatrix@1{\psi_{k} \ar[r]^{H^{1}(\Omega)} &\psi}\right),$$
\textit{with a weighted Poincar\'e inequality.}\\

\noindent
\textsc{Proof.} Let $\left(\psi_{k}\right)_{k \in\mathbb{N}} \in (C^{\infty}_{0})^{\mathbb{N}}$  a Cauchy sequence with respect to the form $a$ whose elements satisfy the boundary condition  $F(\psi_{k})=-\psi_{k}$. Then we have
$$\int_{M} \left|\widehat{\mathfrak{D}}\psi_{k}\right|^{2}_{\gamma} = 
\int_{}\left(\left|\widehat{\nabla}\psi_{k}\right|^{2}_{\gamma}+\left\langle \widehat{\mathfrak{R}}\psi_{k},\psi_{k}\right\rangle_{\gamma}\right)+\frac{1}{2}\int_{\partial M}\left\langle e_{0}\cdot\vec{k}\cdot\psi_{k},\psi_{k}\right\rangle
 \quad ,$$
and thus thanks to the weighted Poincar\'e inequality
$$\forall\Omega\subset M \quad |\Omega|<\infty \quad \xymatrix@1{\psi_{k} \ar[r]^{L^{2}(\Omega)} &\psi} \quad \text{and} \quad \xymatrix@1{\widehat{\nabla}\psi_{k} \ar[r]^{L^{2}(M)}& \rho}.$$
Now let us take a $\varphi\in C^{1}_{0}$ such that $\text{Supp}\varphi\subset K\subset(M \setminus \partial M)$ ($K$ compact without boundary) and then
$$\xymatrix{
 \int_{K}\left\langle \widehat{\nabla}^{*}\varphi,\psi_{k}\right\rangle \ar[d]^{k\rightarrow \infty} \ar@{=}[r] &  \int_{K}\left\langle \varphi,\widehat{\nabla}\psi_{k}\right\rangle \ar[d]^{k\rightarrow \infty} \\
 \int_{K}\left\langle \widehat{\nabla}^{*}\varphi,\psi\right\rangle\ar@{=}[r]&  \int_{K}\left\langle \varphi,\rho\right\rangle },$$
and therefore $\rho =\widehat{\nabla}\psi$ in the distributional sense.\hspace{\stretch{4}}$\square$\\

\noindent
We consider the linear form $l$ on $H_{-}(a)$ defined by
$$l(\psi)= \int_{M}\left\langle\widehat{\mathfrak{D}}(f\mathcal{A}\sigma),\widehat{\mathfrak{D}}\psi\right\rangle_{\gamma}\text{d}\mu_{g}.$$
It still is a continuous linear form on the Hilbert space $H_{-}(a)$ (it is complete since the condition $F(\psi)=-\psi$ is closed) and applying again Lax-Milgram theorem we get the existence of a unique $\xi_{0}\in H(a)$ such that $l=a(-\xi_{0},\cdot)$. In other words
$$\forall\psi\in H_{-}(a) \quad  \int_{M}\left\langle\chi,\widehat{\mathfrak{D}}\psi \right\rangle=0,$$
where we have  set $\xi= f\mathcal{A} \sigma +\xi_{0}$ and $\chi = \widehat{\mathfrak{D}}\xi$.\\

\noindent
\textbf{Remark.} We have for any compactly supported smooth spinor fields $\varphi_{k},\, k=1,2$, the integration by parts formula
$$ \int_{M}\left\langle \varphi_{1},\widehat{\mathfrak{D}}\varphi_{2}\right\rangle=\int_{M}\left\langle\widehat{\mathfrak{D}}^{*}\varphi_{1},\varphi_{2} \right\rangle+\int_{\partial M}\left\langle\nu\cdot \varphi_{1},\varphi_{2} \right\rangle.$$\hspace{\stretch{4}}$\square$\\

\noindent
For any $\psi \in C_{0}^{1}$ we have
$$ \int_{M}\left\langle \chi,\widehat{\mathfrak{D}}\psi\right\rangle=0=\int_{M}\left\langle\widehat{\mathfrak{D}}^{*}\chi,\psi \right\rangle+\int_{\partial M}\left\langle\nu\cdot \chi,\psi \right\rangle.$$
But remembering that $C_{0}^{\infty}(M\setminus\partial M)$ the space of smooth spinor fields compactly supported in $M\setminus \partial M$ is dense in $L^{2}(M)$ then we obtain that $\widehat{\mathfrak{D}}^{*}\chi=0$ and $\chi\in H_{+}(a)=\left\{\psi\in H(a)/ F(\psi)=+\psi\right\}$. By ellipticity $\chi$ is smooth and $\widehat{\mathfrak{D}}^{k}\chi\in L^{2}(M)$ for every $k\in \mathbb{N}$. Finally we notice that
\begin{eqnarray*}
	\int_{M}\left|\widehat{\mathfrak{D}}\chi\right|^{2}&=&\int_{M}\left\langle \widehat{\mathfrak{D}}^{*}\widehat{\mathfrak{D}}\chi,\chi\right\rangle + \int_{\partial M}\left\langle \nu\cdot\widehat{\mathfrak{D}}\chi,\chi\right\rangle\\
	&=& 0 + \int_{\partial M}\left\langle -\textbf{i}n\nu\cdot\chi,\chi\right\rangle\\
	&=& -n\int_{\partial M}\left|\chi\right|^{2},
\end{eqnarray*}
and therefore $\widehat{\mathfrak{D}}\chi=0$ which implies that $\chi=0$. We can conclude with the\\

\noindent
\textbf{Proposition.} {\it  For every $\sigma \in IKS(\Sigma)$ there exists a unique $\xi_{0}\in H_{-}(a)$ such that\\ $\xi=f\mathcal{A}\sigma +\xi_{0} \in \text{Ker} \widehat{\mathfrak{D}}  \cap  H_{-}(a) $ and
\begin{eqnarray*}
\mathcal{H}(V_{\sigma},\alpha_{\sigma})&= &\lim_{r\rightarrow+\infty} \int_{S_{r}}\left\langle \widehat{\nabla}_{\mathcal{A}\nu_{r}}\xi+ \mathcal{A}\nu_{r}\cdot\widehat{\mathfrak{D}}\xi, \xi\right\rangle_{\gamma}\\
&=&  4\int_{M}\left(\left|\widehat{\nabla}\xi\right|^{2}_{\gamma}+\left\langle\widehat{\mathfrak{R}}\xi,\xi\right\rangle_{\gamma}\right)+2\int_{\partial M}\left\langle e_{0}\cdot\vec{k}\cdot\psi,\psi\right\rangle\geq 0.
\end{eqnarray*} } \\

\subsection{Imaginary Killing Spinors}

In general (that is to say whatever the dimension), the space of imaginary Killing spinors of AdS$^{n,1}$ is a finite dimensional complex vector space which trivializes the spinor bundle of AdS$^{n,1}$ (cf. \cite{CGLS} for instance). Consequently there exists an integer $d$ depending upon the dimension $n+1$, such that $IKS(\Sigma)\cong \Bbb C^{d}$ where $\Sigma $ still is the spinor bundle of AdS$^{n,1}$ restricted to a standard hyperbolic slice. Thereby there exists a quadratic Hermitian application
$$\mathcal{K}:\xymatrix{
\Bbb C^{d} \ar[r]^-{\sim}& IKS(\Sigma) \ar[r] &   \Bbb R^{n,1}\oplus \mathfrak{so}(n,1)
}$$
which is difficult to explicite except when the dimension of the slice is $n=3$ (because of exceptional isomorphisms for Lie groups). In order to make out the meaning of the non-negativity of the energy-momentum, we will consider in this section the particular case of dimension $n=3$. \\

\noindent
The aim of this section is to solve explicitely the  Killing equation of section 2.2. As a matter of fact, representation theory provides us good candidates for the imaginary Killing spinors. Thanks to Schur's lemma, we have an isomorphism
$$\begin{array}{ccc}
\mathbb{C}^{2} & \longrightarrow &\text{Hom}^{\text{SU}(2)}(\mathbb{C}^{2},\mathbb{C}^{2}\oplus\mathbb{C}^{2})\\
\binom{z_{1}}{z_{2}}	& \longmapsto & \binom{z_{1}\text{I}_{2}}{z_{2}\text{I}_{2}}
\end{array},$$
We are now considering two families of spinors which are derived from representation theory. To this end, we will denote $w \otimes z \in \mathbb{C}^{2} \otimes \text{Hom}^{\text{SU}(2)}(\mathbb{C}^{2},\mathbb{C}^{2}\oplus\mathbb{C}^{2})$ thanks to the isomorphism above.\\

\noindent
\textbf{Definition.} \textit{Let} $w \otimes z \in \mathbb{C}^{2}\otimes\mathbb{C}^{2}  $ \textit{and set} $ \sigma^{-1}_{w\otimes z}[\tilde{g}]=\left\{\tilde{g}, z(\tilde{g}^{-1}w)\right\}, \ \sigma^{*}_{w\otimes z}[\tilde{g}]=\left\{\tilde{g}, z(\tilde{g}^{*}w)\right\}.$\\

\noindent
Let us consider a spinor field $\tau \in \Gamma(\Sigma)$ and a vector field $X \in \Gamma(\mathbb{T})$ tangent to $\mathbb{H}^{3}$. We can write $ \tau[\tilde{g}] = \left\{\tilde{g}, v(\tilde{g})\right\}$ and $X[g]= \left\{g,\zeta(g)\right\}$, where $v: \mathbb{H}^{3}\longrightarrow  \mathbb{C}^{4}$ and $ \zeta : \mathbb{H}^{3}\longrightarrow \mathfrak{G} $ are respectively $\text{SO}(3)$ and $\text{SU}(2)$-équivariant functions. We can now differentiate $\tau$ in the direction of $X$ and write down
$$( D_{X}\tau )[\tilde{g}]= \{ \tilde{g}, v_{*}(X)_{[\tilde{g}]}+ \tilde{\rho}_{*} \circ s^{*}\theta (\zeta)_{[\tilde{g}]}v[\tilde{g}]\},$$where $\theta$ is the connection 1-form of the bundle of $\text{SL}(2,\mathbb{C})$-frames, restricted to  $\mathbb{H}^{3}$. If one remembers that $\theta$ is only the projection on the first factor in the decomposition $\mathfrak{sl}_{2}(\mathbb{C})\cong \mathfrak{su}(2) \oplus \mathfrak{G}$, we can conclude that $\tilde{\rho}_{*} \circ s^{*}\theta (\zeta)_{[\tilde{g}]}v[\tilde{g}]$ vanishes. Besides we will apply this formula to spinors in $\left\{\sigma^{-1}_{w\otimes z},\sigma^{*}_{u\otimes z}, w,u \in \mathbb{C}^{2} \right\}$ so that we can only derive at the point $\tilde{g}=1$ unity in $\text{SL}(2,\mathbb{C})$ since we have the\\

\noindent
\textbf{Proposition.} \textit{The set} $\left\{\sigma^{-1}_{w\otimes z},\sigma^{*}_{u\otimes z}, w,u \in \mathbb{C}^{2} \right\}$ \textit{is stable under the} $\text{SL}(2,\mathbb{C})$ \textit{action. More precisely for every} $\tilde{e}\in \text{SL}(2,\mathbb{C})$ \textit{we have} $\tilde{e}* \sigma^{-1}_{w\otimes z} = \sigma^{-1}_{\tilde{e}w\otimes z}$ \textit{and} $\tilde{e}*\sigma^{*}_{u\otimes z} = \sigma^{*}_{(\tilde{e}^{*})^{-1}u\otimes z}.$\\

\noindent
We obtain
$$\left\{\begin{array}{lll}
	( D_{X}\sigma^{-1}_{w\otimes z} )[1]&=&\left\{1, -z( \zeta w)\right\}\\
	( D_{X}\sigma^{*}_{u\otimes z} )[1]&=&\left\{1, z( \zeta u)\right\}
\end{array}\right. ,$$
where $\zeta=\zeta(1)$. We also compute the Clifford action of X on $\sigma^{-1}_{w\otimes z},\sigma^{*}_{u\otimes z}$ at the point 1:
$$\left\{\begin{array}{lll}
X \cdot \sigma^{-1}_{w\otimes z}[1]&=&\left\{1, \Theta(\zeta)z( w)\right\}\\
X \cdot \sigma^{*}_{u\otimes z}[1]&=&\left\{1, \Theta(\zeta)z( u)\right\}
\end{array}\right. .$$
We must precise $ \Theta_{|\mathfrak{G}}: \zeta  \longmapsto \left(\begin{array}{cc}
0 & 2\zeta \\
-2 \zeta & 0	
\end{array}\right)$, and if we introduce the sections $\sigma^{-1}_{w\otimes \binom{1}{-i}}$ and $\sigma^{*}_{w\otimes \binom{1}{i}}$, for any $w \in \mathbb{C}^{2}$, we have on one hand
$$\left\{\begin{array}{lllll}
-\frac{i}{2}X \cdot \sigma^{-1}_{w\otimes\binom{1}{-i} }[1]&=&-i\left\{1, -i \zeta w\oplus -\zeta w \right\} & =&\left\{1, - \zeta w\oplus i\zeta w \right\} \\
-\frac{i}{2}X \cdot \sigma^{*}_{u\otimes \binom{1}{i}}[1]&=&-i\left\{1, i \zeta u\oplus -\zeta u\right\} &=&\left\{1, \zeta u\oplus i\zeta u \right\}
\end{array}\right. ,$$
and on the other hand
$$\left\{\begin{array}{cll}
	\left( D_{X}\sigma^{-1}_{w\otimes \binom{1}{-i}} \right)[1]&=&\left\{1, - \zeta w \oplus i \zeta w \right\}\\
	\left( D_{X}\sigma^{*}_{u\otimes \binom{1}{i}} \right)[1]&=&\left\{1, \zeta u \oplus i \zeta u \right\}
\end{array}\right. .$$
 Since $ \left\{\sigma^{-1}_{w\otimes \binom{1}{-i}}+\sigma^{*}_{u\otimes \binom{1}{i}} /  w, u\in \mathbb{C}^{2} \right\}$ is a 4-dimensional complex vector space, we obviously obtain the\\ 

\noindent
\textbf{Proposition.} \textit{The space of imaginary Killing spinors denoted by} $IKS(\Sigma)$ \textit{is generated by}
$$\left\{\sigma^{-1}_{w\otimes \binom{1}{-i}},\sigma^{*}_{u\otimes \binom{1}{i}}, w,u\in \mathbb{C}^{2} \right\}.$$
Let $\sigma$ an imaginary Killing spinor and set $V_{\sigma}:= <\sigma,\sigma>$ which is a function on  $\mathbb{H}^{3}$, and if $e_{0}$ denotes a unit normal of $\mathbb{H}^{3}$ in AdS, we set $\alpha_{\sigma}(Y):= \left\langle Y \cdot e_{0} \cdot\sigma,\sigma \right\rangle$ which is a real 1-form on $\mathbb{H}^{3}$. The goal of the two next paragraphs is  to define some $\text{SL}(2,\Bbb C)$-equivariant application
$$\begin{array}{lrll}
 \mathcal{K}:& IKS(\Sigma)\cong \mathbb{C}^{2}\oplus \mathbb{C}^{2} & \longrightarrow & \left(\mathfrak{M}\oplus\mathfrak{sl}_{2}(\mathbb{C})\right)^{*\mathbb{R}} \\
& w\oplus u  & \longmapsto &  \mathcal{K}_{w\oplus u }:= (V_{w\oplus u }\oplus\alpha_{w\oplus u }).
\end{array}$$

\begin{center}
\textbf{\large{The functions $V_{\sigma}$}}
\end{center}

\noindent
We compute the functions $V_{\sigma}$ which are by definition
$$
\begin{array}{lll}
V_{\sigma}[\tilde{g}] &=& 	\left|\sigma^{-1}_{w\otimes \binom{1}{-i}}[\tilde{g}]\right|^{2}_{\mathbb{C}^{4}}+ \left|\sigma^{*}_{u\otimes \binom{1}{i}}[\tilde{g}]\right|^{2}_{\mathbb{C}^{4}}+
2 \Re e\left(\left\langle \sigma^{-1}_{w\otimes \binom{1}{-i}}[\tilde{g}],\sigma^{*}_{u\otimes \binom{1}{i}}[\tilde{g}]\right\rangle_{\mathbb{C}^{4}}\right)\\
&=& 2\left|\tilde{g}^{-1}w\right|^{2}_{\mathbb{C}^{2}} + 2\left|\tilde{g}^{*}u\right|^{2}_{\mathbb{C}^{2}}
\end{array}$$

\noindent
\textbf{Remark.} $\sigma^{-1}_{w\otimes \binom{1}{-i}}$ and $\sigma^{*}_{u\otimes \binom{1}{i}}$ are orthogonal spinors for every $ u,w \in\mathbb{C}^{2}$.\hspace{\stretch{4}}$\square$\\

\noindent
If $\tilde{g}\in \text{SL}(2,\mathbb{C})$, the corresponding base point is  $\tilde{g}\tilde{g}^{*} \in \mathbb{H}^{3}\subset \mathfrak{M} \cong \mathbb{R}^{3,1}$ whose coordinates are given by $ (x_{k})_{k=0}^{3}= \Lambda^{-1}(\tilde{g}\tilde{g}^{*})$.\\

\noindent
\textbf{Proposition.}  $V_{\sigma}$ \textit{is a causal element of} $N_{b}.$\\

\noindent
\textsc{Proof.} Let $U=\binom{u_{1}}{-\overline{w_{2}}}\in \mathbb{C}^{2}, V=\binom{u_{2}}{\overline{w_{1}}} \in\mathbb{C}^{2}$. We notice that
$$V_{\sigma}[\tilde{g}]  =  x_{0} (|U|^{2}+|V|^{2}) +x_{1}(|U|^{2}-|V|^{2})  + 2x_{2}\Re e (<U,V>)- 2x_{3}\Im m (<U,V>),$$
so that the norm of $V_{\sigma}$ is
 $\left|V_{\sigma}[\tilde{g}]\right|^{2}= 4\left( |<U,V>|^{2}-|U|^{2}|V|^{2}\right)\leq 0$, thanks to the Cauchy-Schwarz inequality for the standard Hermitian form on $\mathbb{C}^{2}$.\hspace{\stretch{4}}$\square$\\
 
\noindent
More conceptually we see that $V_{\sigma}[\tilde{g}] = 2 ( w^{*}\widehat{W}w+ u^{*}Wu)$, where we have set ${W{:=} \tilde{g}\tilde{g}^{*}{\in} \mathbb{H}^{3} {\subset} \mathfrak{M}}$. Thereby we can define by extension an application 
$$\begin{array}{lll}
\mathbb{C}^{2}\oplus \mathbb{C}^{2}& \longrightarrow &\mathfrak{M}^{*} \\
w \oplus u & \longmapsto & \left( V_{w \oplus u}: W \mapsto  2 ( w^{*}\widehat{W}w+ u^{*}Wu) \right)
\end{array}.$$

\begin{center}
\large{\textbf{The 1-forms $\alpha_{\sigma}$}}
\end{center}
\noindent
The positively oriented unit normal of $\mathbb{H}^{3}$ in AdS is  given by $e_{0}[\tilde{g}]=\left\{ \tilde{g}, \frac{1}{2} \mu(\tilde{g})I_{2}\right\}$ and for any $\xi\in \mathfrak{G}$ satisfying $-\text{det}\xi=1$ we set $ X^{\xi}[\tilde{g}]=\left\{ \tilde{g}, \frac{1}{2} \mu(\tilde{g})\xi \right\}$. Just remember that $\alpha_{\sigma}(X^{\xi})_{[\tilde{g}]}:= \left\langle X^{\xi} \cdot e_{0} \cdot\sigma,\sigma \right\rangle_{[\tilde{g}]}$. As we suppose that $\sigma\in IKS(\Sigma)$, we can easily compute the first derivative of $\alpha_{\sigma}$
$$ D_{X^{\eta}} \alpha_{\sigma}(X^{\xi})_{[\tilde{g}]}= \frac{\textbf{i}}{2}\left\langle ( X^{\eta}\cdot X^{\xi}-X^{\xi}\cdot X^{\eta})\cdot e_{0} \cdot \sigma, \sigma \right\rangle_{[\tilde{g}]},$$
which is a real skew symmetric 2-form and hence $\alpha_{\sigma}$ is a Killing form on $\mathbb{H}^{3}$ . From now on we set $ \alpha_{\sigma}=(\alpha_{\sigma})_{1}$ and $ D\alpha_{\sigma}=(D \alpha_{\sigma})_{1},$ 
that we will write as function of $w\oplus u$. After some computations we find

$$\left\{\begin{array}{lll}
\alpha_{\sigma}(\xi) & = & 2(w^{*}\xi u + u^{*} \xi w)  \\
D \alpha_{\sigma}(\eta,\xi) & = & (w^{*}(\xi\eta-\eta\xi)  u - u^{*}(\xi\eta-\eta\xi)  w)
\end{array} \right. .$$
We have to notice  that $\xi\eta-\eta\xi \in \textbf{i}\mathfrak{G}$ so that $D \alpha_{\sigma}$ is naturally a linear form on $\textbf{i}\mathfrak{G}$. As a consequence we define, thanks to the Killing 1-form $\alpha_{\sigma}$, the following application
$$\begin{array}{cll}
\mathbb{C}^{2}\oplus \mathbb{C}^{2} & \longrightarrow & \mathfrak{sl}_{2}(\mathbb{C})^{*\mathbb{R}} \\
w\oplus u  & \longmapsto &  \left(\alpha_{w\oplus u }: \xi \mapsto 2(w^{*}\xi u + u^{*}\xi^{*} w)\right),
\end{array}$$\\
where $^{*\mathbb{R}}$ stands for the duality with respect to the reals. We then define 
$$\mathcal{K}_{w\oplus  u}=V_{w\oplus u }\oplus  \alpha_{w\oplus u }$$ and conclude with the\\

\noindent
\textbf{Proposition.} \textit{The application} $\mathcal{K}$ \textit{is} $\text{SL}(2,\mathbb{C})$\textit{-equivariant. More precisely, for every} ${\tilde{e}\in\text{SL}(2,\mathbb{C})}$
$$\mathcal{K}_{\tilde{e}*(w\oplus u)}= \left(V_{w\oplus u }\circ \mu(\tilde{e}^{-1})\right) \oplus \left(\alpha_{w\oplus u }\circ \text{Ad}(\tilde{e}^{*})\right).$$\\

\noindent
\textsc{Proof.} We must compute for every $W\in \mathfrak{M}$ and $\xi \in \mathfrak{sl}_{2}(\mathbb{C})$
$$
\begin{array}{lll}
\mathcal{K}_{\tilde{e}*(w\oplus u)}(W,\xi) &=& \mathcal{K}_{\tilde{e}w\oplus (\tilde{e}^{*})^{-1}u}(W,\xi) \\
&=& 2(w^{*}\tilde{e}^{*}\widehat{W}\tilde{e}w + u^{*}\tilde{e}^{-1}W(\tilde{e}^{*})^{-1}u + w^{*}\tilde{e}^{*}\xi (\tilde{e}^{*})^{-1}u + u^{*}\tilde{e}^{-1}\xi^{*} \tilde{e}w )\\
&=& V_{w\oplus u }\circ \mu(\tilde{e}^{-1})(W)\oplus  \alpha_{w\oplus u }\circ \text{Ad}(\tilde{e}^{*})(\xi).
\end{array}$$ \hspace{\stretch{4}}$\square$

\noindent
\textbf{Remark. The norm of imaginary Killing spinors} \\
Classical considerations on Lie algebras show that $\mathfrak{so}(3,2) $ endowed with its Killing form, is isometric to  $ (\mathfrak{M},- \text{det}) \oplus (\mathfrak{sl}_{2}(\mathbb{C}), -\Re e(\text{det}))$ which is a 10-dimensional real vector space of signature (6,4). The norm of $\mathcal{K}(w\oplus u)$ with respect to the Killing form is, up to a multiplicative and positive constant $\left|\mathcal{K}(w\oplus u) \right|^{2}= |<U,V>|^{2}-|U|^{2}|V|^{2}+\Re e ( \chi^{2})$, where we have set $\chi=\overline{u_{1}}w_{1} + \overline{u_{2}}w_{2}$. Besides, if $V_{w\oplus u}$ is isotropic in $\mathfrak{M}$ then $\alpha_{w\oplus u}$ and $\mathcal{K}(w\oplus u)$ are also isotropic respectively in $\mathfrak{sl}_{2}(\mathbb{C})^{*}$ and $\left(\mathfrak{M} \oplus \mathfrak{sl}_{2}(\mathbb{C})\right)^{*}$.  Indeed the equality case in the Cauchy-Schwarz inequality occurs if and only if \textit{U} and \textit{V} satisfy $\text{det}_{\mathbb{C}^{2}}(U,V)=\bar{\chi}=0.$\hspace{\stretch{4}}$\square$\\

\subsection{End of the Proof}

Whatever the dimension is, we obtain a Hermitian application 
$$Q :\xymatrix{
\Bbb C^{d} \ar[r]^-{\mathcal{K}} & \Bbb R^{n,1}\oplus \mathfrak{so}(n,1)  \ar[r]^-{\mathcal{H}} & \Bbb R
} ,$$
which has to be non-negative in vertue of the non-negativity results of sections 3.2.1 and 3.2.2. This completes the proof of the positivity theorem stated in section 1.3.\\
In dimension $n=3$, we can be more specific giving the explicite formula of $Q$ in terms of the components of the energy-momentum $\mathcal{H}$. More precisely, on one hand we have found a quadratic application
$$\begin{array}{lrll}
 \mathcal{K}:& IKS(\Sigma)\cong \mathbb{C}^{2}\oplus \mathbb{C}^{2} & \longrightarrow & \left(\mathfrak{M}\oplus\mathfrak{sl}_{2}(\mathbb{C})\right)^{*\mathbb{R}}\cong \text{Ker d}\Phi^{*}_{(b,0)} \\
& w\oplus u  & \longmapsto &  (V_{w\oplus u }\oplus\alpha_{w\oplus u })
\end{array},$$
which is $\text{SL}(2,\mathbb{C})$-equivariant. On the other hand we know that the energy-momentum functional $\mathcal{H}$ can be seen as a real linear form on $\left(\mathfrak{M}\oplus\mathfrak{sl}_{2}(\mathbb{C})\right)^{*\mathbb{R}}$ that is to say, as a vector   $\mathcal{H}=M\oplus\Xi\in\mathfrak{M}\oplus\mathfrak{sl}_{2}(\mathbb{C})$. In the following, we will adopt the notations $\Xi=N\oplus \textbf{i}R \in \mathfrak{G}\oplus\textbf{i}\mathfrak{G}$, and  $M=\Lambda (m_{0},m), \,  N= \Lambda(0,n),\, R=\Lambda(0,r)$, where $\Lambda$ is the isomorphism defined in section 2.4. Now applying the non-negativity result of section 3.2.1 or 3.2.2, we know that (even if our AdS-asymptotically hyperbolic manifold has a compact boundary such that $\vec{k}$ is causal and positively oriented)
$$ \forall \sigma \in IKS(\Sigma ) \qquad \mathcal{H}(V_{\sigma },\alpha _{\sigma }) \geq  0.$$
In other words, for each $w\oplus u \in \Bbb C^{4}$, we have  $\mathcal{H}(\mathcal{K}_{w\oplus u })\geq 0$. But the complete study of $IKS(\Sigma )$ of section 3.3 implies that actually
\begin{eqnarray*}
\mathcal{H}(\mathcal{K}_{w\oplus u })&=& V_{w\oplus u }(M) +\alpha_{w\oplus u }(\Xi)\\
&=& 2( w^{*}\widehat{M}w+ u^{*}Mu) + 2(w^{*}\Xi u + u^{*}\Xi^{*} w),
\end{eqnarray*}\\
and consequently the application $ w\oplus u \longmapsto \mathcal{H}(\mathcal{K}_{w\oplus u })$ is a Hermitian form on  $\mathbb{C}^{2}\oplus \mathbb{C}^{2}$ whose matrix is
$$Q=2
\left(\begin{array}{cc}
	\widehat{M}& \Xi\\
	\Xi^{*} & M
\end{array}\right)=2
\left(\begin{array}{cc}
\Lambda (m_{0},-m) & \Lambda(0,n)+ \textbf{i}\Lambda(0,r) \\
\Lambda(0,n)- \textbf{i}\Lambda(0,r) & \Lambda (m_{0},m) 
\end{array} \right).$$
It is easy to conclude since we have the identity
$$ \forall w\oplus u\in \Bbb C^{4} \quad \mathcal{H}(V_{w\oplus u }\oplus \alpha _{w\oplus u })=Q(w\oplus u,w\oplus u) \geq  0,$$
which ends  the proof of the \\

\noindent
\textbf{Positive Energy-Momentum Theorem.} \textit{Let} $(M^{n},g,k)$ \textit{be an AdS-asymptotically hyperbolic spin Riemannian manifold satisfying the decay conditions stated in section 1.2 and the following conditions\\
(i)} $\left\langle (f, \alpha),(\Phi(g,k)- \Phi(b,0))\right\rangle \in L^{1}(M,\text{dVol}_{b})$ \textit{for every} $(f,\alpha)\in N_{b}\oplus \mathfrak{Kill}(M,b),$ \\
\textit{(ii) the relative version of the dominant energy condition (cf. section 2.2) holds, that is to say} $(\Phi(g,k)- \Phi(b,0))$ \textit{is a positively oriented causal (n+1)-vector along $M$,\\
(iii) in the case where M has a compact boundary $\partial M$, we assume moreover that $\vec{k}$ is causal and  positively oriented along $\partial M$.\\
Then there exists  a (hardly explicitable) map} $\Bbb R^{n,1} \oplus \mathfrak{so}(n,1) \longrightarrow  \text{Herm}(C^d)$ {\it which sends, under the assumptions (i-iii), the energy-momentum on a non-negative Hermitian form Q.\\
Moreover, when n=3, we can explicite Q in terms of the components of the energy-momentum as described above.}\\

\noindent
The end of this section is devoted to the 3-dimensional case.\\
As the invariance of the non-negativity of $Q$ under asymptotic hyperbolic isometries, was proved in \cite{ChN}, one can be interested in the description of the orbit of the energy-momentum under the action of $\text{SL}(2,\mathbb{C})$.\\

\noindent
\textbf{Proposition.} {\it If M is timelike, there exists a (non-unique) representative element of the orbit of $\mathcal{H}=M\oplus \Xi$ under the natural action (cf. section 3.3) of $\text{SL}(2,\mathbb{C})$ on $\mathfrak{M}\oplus\mathfrak{sl}_{2}(\mathbb{C})$ which can be written
$$m_{0} \left(
\begin{array}{cc}
1 & 0\\
0 & 1
\end{array}
\right)  \oplus n_{1}  \left(
\begin{array}{cc}
1 &0 \\
0 & -1
\end{array}
\right)  \oplus \textbf{i}
\left(
\begin{array}{cc}
r_{1} & r_{2} \\
r_{2} & -r_{1}
\end{array}
\right), \ m_{0}, n_{1}, r_{1}, r_{2}\in \mathbb{R}.$$
The positive energy-momentum theorem then reduces to $m_{0} \geq \sqrt{(|n_{1}|+|r_{2}|)^{2}+r_{1}^{2}}.$}\\

\noindent
\textsc{Proof.} Let us suppose that $M\in\mathfrak{M}$ is timelike. Thus considering the action of  $\text{SL}(2,\mathbb{C})$ on $\mathfrak{M}\oplus\mathfrak{sl}_{2}(\mathbb{C}) $ (cf. section 3.3), then there exists an element in the orbit of $\mathcal{H}$ that can be written $ m_{0} \left(
\begin{array}{cc}
1 & 0\\
0 & 1
\end{array}\right)\oplus \Xi'$. Since the isotropy group of  $\left(
\begin{array}{cc}
1 & 0\\
0 & 1
\end{array}\right)$ is $\text{SU}(2)$ whose action on $\mathfrak{G}$ is transitive, then there exists an element in the orbit of $\mathcal{H}$ that can be written $ m_{0} \left(
\begin{array}{cc}
1 & 0\\
0 & 1
\end{array}\right)\oplus  n_{1}  \left(
\begin{array}{cc}
1 &0 \\
0 & -1
\end{array}
\right)  \oplus \textbf{i} R'$. But the isotropy group of  $\left(
\begin{array}{cc}
1 &0 \\
0 & -1
\end{array}
\right)$ is the one parameter group $ \left\{   \left(
\begin{array}{cc}
e^{\textbf{i}\theta} & 0  \\
0 & e^{-\textbf{i}\theta}
\end{array}
\right),\theta\in\mathbb{R}    \right\}$. Finally there exists an element (not unique since the isotropy group of $\left(
\begin{array}{cc}
0 & 1 \\
1 & 0
\end{array}
\right)$ is isomorphic to $\mathbb{Z}_{2}$) in the orbit of $\mathcal{H}$ that can be written as announced in the proposition. The corresponding Hermitian matrix  is 
$$Q=2\left(\begin{array}{cccc}
m_{0} & 0 & n_{1}+ \textbf{i}r_{1} & \textbf{i}r_{2} \\
0 & m_{0} &  \textbf{i}r_{2} & - n_{1}- \textbf{i}r_{1}\\ 
 n_{1}- \textbf{i}r_{1} & - \textbf{i}r_{2}& m_{0} & 0 \\
- \textbf{i}r_{2} & - n_{1}+\textbf{i}r_{1} & 0 & m_{0}
\end{array}\right).$$
Since $Q$ is non negative we have
\begin{eqnarray*}
m_{0}& \geq& 0 \\
m_{0}(m_{0}^{2}-(n_{1}^{2}+r_{1}^{2}+r_{2}^{2})) &\geq& 0\\
(m_{0}^{2}-(n_{1}^{2}+r_{1}^{2}+r_{2}^{2}))^{2} &\geq& 4(n_{1}r_{2})^{2},
\end{eqnarray*}
which can be summerized with $m_{0} \geq \sqrt{(|n_{1}|+|r_{2}|)^{2}+r_{1}^{2}}$.\hspace{\stretch{4}}$\square$\\

\noindent
{\bf Remark.}
The $\left\{ t=0\right\}$ slices of the Kerr AdS metrics are AdS-asymptotically hyperbolic and parametrized by 2 real parameters: {\it the mass and the angular momentum}. The proposition above then shows that there exists some energy-momenta that could not be obtained by the action of $\text{SL}(2,\mathbb{C})$ on a Kerr-AdS solution. As a consequence, an interesting question would be to find some (new) AdS-asymptotically hyperbolic metrics which have an energy-momentum of the form given in the proposition above with non-zero coefficients $m_{0}, n_{1}, r_{1}, r_{2}$, and which satisfy the dominant energy condition or the (stronger) cosmological vacuum constraints.

\section{Rigidity Theorems}

\noindent
\textbf{Theorem.} {\it Under the assumptions of the positive energy-momentum theorem,}  $\text{tr}Q=0$ {\it implies that  $(M,g,k)$ is  isometrically embeddable in} AdS$^{n,1}$.\\

\noindent
{\sc Proof.}
The vanishing of $\text{tr}Q$ implies from the non-negativity of $Q$, that $Q=0$. Consequently our spinor bundle $\Sigma $ is trivialized by a basis of {\it $\gamma $-imaginary Killing spinors}. We denote by $\xi $ any {\it $\gamma $-imaginary Killing spinors} of this basis. We will need the following spinorial Gauss-Codazzi equation.\\

\noindent
\textbf{Proposition.} \textit{For every} $X ,Y \in \Gamma(TM)$ \textit{we have}
$$ R^{\gamma}_{X,Y}= R^{g}_{X,Y} -\frac{1}{2}\left( \text{d}^{\overline{\nabla}}k(X,Y)\cdot e_{0} + \frac{1}{2}\big(k(X)\cdot k(Y)-k(Y)\cdot k(X)\big) \right)\cdot \quad ,$$
\textit{where} $\cdot$ \textit{denotes the Clifford action with respect to the metric} $\gamma$.\\

\noindent
\textsc{Proof of the proposition.} It is a straightforward computation where we use vector fields $X,Y$ satisfying at the point $\overline{\nabla}_{X}Y=\overline{\nabla}_{Y}X=0$.
\begin{eqnarray*}
\nabla_{X} \nabla_{Y}  & = & \nabla_{X} \left( \overline{\nabla}_{Y} -\frac{1}{2}k(Y)\cdot e_{0}\cdot   \right)\\
&=&  \overline{\nabla}_{X} \overline{\nabla}_{Y}   -\frac{1}{2}k(X)\cdot e_{0}\cdot\overline{\nabla}_{Y}\\
&&-\frac{1}{2}\left( \nabla_{X}k(Y)\cdot e_{0}\cdot + k(Y)\cdot (\nabla_{X}e_{0})\cdot +k(Y)\cdot e_{0}\cdot \nabla_{X} \right)\\
&=& \overline{\nabla}_{X} \overline{\nabla}_{Y}  -\frac{1}{2}\left( k(X)\cdot e_{0}\cdot\overline{\nabla}_{Y}+k(Y)\cdot e_{0}\cdot\overline{\nabla}_{X} -k \circ k (X,Y)e_{0} \cdot \right)\\
&& -\left(\frac{1}{2} \overline{\nabla}_{X}k(Y)\cdot e_{0}  -\frac{1}{4} k(Y)\cdot k(X) \right)\cdot \quad ,
\end{eqnarray*}
and the curvature formula above follows.\hspace{\stretch{4}}$\square$\\

\noindent
Using the fact that $\xi $ is a $\gamma$-imaginary Killing spinor one gets 
$$ \left\langle R^{g}_{X,Y}\xi -\frac{1}{4}\bigg( X\cdot Y -Y\cdot X + k(X)\cdot k(Y)-k(Y)\cdot k(X) \bigg)\cdot\xi,\xi  \right\rangle  =\frac{1}{2}\left\langle \text{d}^{\overline{\nabla}}k(X,Y)\cdot e_{0} \cdot \xi,\xi  \right\rangle ,$$
where $ \left\langle R^{g}_{X,Y}\xi,\xi  \right\rangle $ and $ \left\langle \bigg( X\cdot Y -Y\cdot X +k(X)\cdot k(Y)-k(Y)\cdot k(X)\bigg)\cdot\xi,\xi  \right\rangle $ are purely imaginary terms whereas $ \left\langle \text{d}^{\overline{\nabla}}k(X,Y)\cdot e_{0} \cdot\xi,\xi  \right\rangle $ is real. As a consequence ${\left\langle \text{d}^{\overline{\nabla}}k(X,Y)\cdot e_{0} \cdot\xi,\xi  \right\rangle=0}$ for any  $\xi$ of our $\gamma$-imaginary Killing spinor basis and so $\text{d}^{\overline{\nabla}}k=0$. This implies  
$$ R^{g}_{X,Y}= \frac{1}{4}\big( X\cdot Y - Y\cdot X + k(X)\cdot k(Y)- k(Y)\cdot k(X)\big) \cdot \quad , $$
and using the natural isomorphism between $\text{C}\ell_{0}(\mathbb{R}^{3,1})$ and $\Lambda ^{2} (\mathbb{R}^{3,1})$ (cf. \cite{LMi} proposition 6.2) we get that 
\begin{eqnarray*}
 R^{g}&=& \frac{1}{2}\big( g\owedge g + k\owedge k \big)\\
\text{d}^{\overline{\nabla}}k&=&0 .
\end{eqnarray*}
Let us denote by $V$ the function $<\xi,\xi>$, $\alpha$ the real 1-form defined by $\alpha(Y)= \left\langle Y \cdot e_{0} \cdot \xi,\xi \right\rangle$. Then the couple $(V,Y):=(V,-\alpha^{\sharp})$ is a Killing Initial Data (KID) \cite{BCh}. If we consider $(\widetilde{M},\tilde{g},\tilde{k})$ the universal Riemannian covering of $(M,g,k)$, then we can make the Killing development of $(\widetilde{M},\tilde{g},\tilde{k})$ with respect to the KID  $(\widetilde{V},\widetilde{Y})$ which by definition is $\mathbb{R}\times \widetilde{M}$ endowed with the Lorentzian metric $\tilde{\gamma}=\left(-\widetilde{N}^{2}+ |\widetilde{Y}|^{2} \right)\text{d}u^{2}+ 2\widetilde{Y}^{\flat}\odot \text{d}u + \tilde{g}$. By construction, $\widetilde{M}$ is embedded in $(\mathbb{R}\times \widetilde{M},\tilde{\gamma})$ with induced metric $\tilde{g}$ and second fundamental form $\tilde{k}$. Besides $\mathbb{R}\times \widetilde{M}$ is the universal covering of $N$, and $\tilde{\gamma} $ which has sectional curvature -1, is a stationary solution of the vacuum Einstein equations with cosmological constant that is to say $G^{\tilde{\gamma}}=\frac{n(n-1)}{2}\tilde{\gamma}$. But $(\widetilde{M},\tilde{g})$ is complete since $(M,g)$ is complete and therefore \cite{And1}  $(\mathbb{R}\times \widetilde{M},\tilde{\gamma})$ is geodesically complete. It follows that  $(\mathbb{R}\times \widetilde{M},\tilde{\gamma})$ is AdS$^{n,1}$ (in vertue of Proposition 23 [p.227] of  \cite{O}). It only remains to show that $M$ is simply connected. We know that $\mathbb{R}\times \widetilde{M}\cong \mathbb{R}^{n+1}$ and thereby using the following compactly supported de Rham cohomology isomorphisms  $\left\{0\right\}=H^{2}_{dR,c}(\mathbb{R}\times \widetilde{M})=H^{2}_{dR,c}(\mathbb{R}^{n+1})=H^{1}_{dR,c}( \widetilde{M})$ (cf. Proposition 4.7 and Corollary 4.7.1 [p. 39] of \cite{BT} for instance), we obtain that $\widetilde{M}$ has only one asymptotic end. This last fact compels the universal covering map $\widetilde{M}\rightarrow  M$ to be trivial and as a consequence $(M,g,k)\equiv (\widetilde{M},\tilde{g},\tilde{k})$ is isometrically embedded in AdS$^{n,1} \equiv (\widetilde{N},\tilde{\gamma})$. This completes the proof of the theorem. \hspace{\stretch{4}}$\square$\\

\noindent
Then end of this section is devoted to weaken the condition defining the rigidity case in dimension $n=3$. Namely we prove the following\\

\noindent
{\bf Theorem.} {\it Let us suppose that $(M^{3},g,k)$ satisfies the assumptions of the positive energy-momentum theorem and that the matrix Q is degenerate. Then there exists some $\widehat{\nabla}$-parallel spinor field $\xi $ such that  $\left\langle \widehat{\mathfrak{R}}\xi,\xi\right\rangle =0$ and consequently $(M,g,k)$ is isometrically embeddable in a stationary pp-wave space-time.\\
If furthermore the constant function $(\xi ,\xi )$ is non-zero then $(M,g,k)$ admits a vacuum Killing development which is a solution of the Einstein equations (with  the cosmological constant -3).}\\

\noindent
{\bf Remark.} A pp-wave space-time is a Lorentzian manifolds such that its stress-energy tensor satisfies $T_{\mu \nu }=\lambda Z_{\mu }\otimes Z_{\nu }$ where $Z^{\mu }$ is an isotropic Killing vector field and $\lambda $ a function on the manifold. Some results were also proved by Siklos in \cite{S} and by Leitner in \cite{L} for Lorentzian manifolds admitting a Killing spinor.\hspace{\stretch{4}}$\square$\\

\noindent
{\sc Proof.} The degenerate character of $Q$ implies the existence of a non-zero $w\oplus u \in \Bbb C^{4}$ and a unique $\xi_{0}$ such that $\xi=f\mathcal{A}\sigma_{w\oplus u}+\xi_{0}$ satisfies the conditions
\begin{eqnarray*}
	\widehat{\nabla}\xi &=& 0\\
	\left\langle \widehat{\mathfrak{R}}\xi,\xi\right\rangle &=&0 .
\end{eqnarray*} 
By the same argument as above we get that ${\left\langle \text{d}^{\overline{\nabla}}k(X,Y)\cdot e_{0} \cdot\xi,\xi  \right\rangle=0}$ (which can also be thought as $\overline{\nabla}_{X}k(Y,\alpha)=\overline{\nabla}_{Y}k(X,\alpha)$). 
Now since $\xi $ is $\widehat{\nabla}$-parallel we get
\begin{eqnarray*}
\Re e \left\langle \sum_{k=1}^{3} e_{k} \cdot R_{X,e_{k}}^{\gamma } \xi, Y\cdot\xi  \right\rangle &=& \frac{1}{4}\Re e \left\langle \sum_{k=1}^{3} e_{k} \cdot (X \cdot  e_{k} - e_{k}\cdot X)\cdot \xi, Y\cdot\xi  \right\rangle \\
&=&  \left\langle X ,Y  \right\rangle \qquad \forall X,Y \in \Gamma(TM).
\end{eqnarray*}
On the other hand a direct computation leads to
\begin{eqnarray*}
\sum_{k=1}^{3} e_{k} \cdot R_{X,e_{k}}^{\gamma } &=&\frac{1}{2} \sum_{l,m=1}^{3} \text{d}^{\overline{\nabla}}k(X,e_{l},e_{m}) e_{l}\cdot e_{0}\cdot e_{m}\cdot  + \frac{1}{2}\sum_{l,m=1}^{3} R^{\gamma}(X,e_{l},e_{l},e_{m})e_{l}\cdot e_{l}\cdot e_{m}\cdot\\
&=& \frac{1}{2} \sum_{l,m=1}^{3} \text{d}^{\overline{\nabla}}k(X,e_{l},e_{m}) e_{l}\cdot e_{0}\cdot e_{m}\cdot\\
&&- \frac{1}{2}\sum_{l,m=1}^{3} \left\{ R^{g}(X,e_{l},e_{l},e_{m}) -k(X,e_{l})k(e_{l},e_{m}) +  k(X,e_{m})k(e_{l},e_{l}) \right\} \cdot e_{m}\cdot\\
&=& \frac{1}{2} \sum_{l,m=1}^{3} \text{d}^{\overline{\nabla}}k(X,e_{l},e_{m}) e_{l}\cdot e_{0}\cdot e_{m}\cdot - \frac{1}{2}(E(X)-2X)\cdot \quad,
\end{eqnarray*}
where we have set $E=\text{Ric}^{g}+2g +(\text{tr}k)k -k\circ k$. It is then clear that
$$  \Re e \left\langle 	\sum_{k,l=1}^{3} \text{d}^{\overline{\nabla}}k(X,e_{l},e_{m}) e_{l}\cdot e_{0}\cdot e_{m} \cdot\psi,Y\cdot\psi \right\rangle=VE(X,Y) $$
In the following computation we will set $Y=e_{s}$. We recall that
\begin{eqnarray*}
	\Re e \left\langle 	\sum_{l,m=1}^{3} \text{d}^{\overline{\nabla}}k(X,e_{l},e_{m}) e_{l}\cdot e_{0}\cdot e_{m} \cdot\psi,e_{s}\cdot\psi \right\rangle &=& \sum^{3}_{l,m=1}\text{d}^{\overline{\nabla}}k(X,e_{l},e_{m})\Re e\left\langle  e_{l}\cdot e_{0}\cdot e_{m}\cdot\psi,e_{s}\cdot \psi \right\rangle\\
	&=& \sum^{3}_{l=1}\text{d}^{\overline{\nabla}}k(X,e_{l},e_{l})\Re e\left\langle  e_{l}\cdot e_{0}\cdot e_{l}\cdot\psi,e_{s}\cdot \psi \right\rangle\\
	&& + \sum_{l\neq m}\text{d}^{\overline{\nabla}}k(X,e_{l},e_{m})\Re e\left\langle  e_{l}\cdot e_{0}\cdot e_{m}\cdot\psi,e_{s}\cdot \psi \right\rangle\\
	&=& (I + II)(X,e_{s}),
\end{eqnarray*}
and we will treat $I$ and $II$ seperately for convenience. The easiest one is 

\begin{eqnarray*}
	I (X,e_{s})&=&-\left\langle e_{s}\cdot e_{0}\cdot\psi,\cdot \psi \right\rangle\sum^{3}_{l=1}\bigg(\overline{\nabla}_{X}k(e_{l},e_{l})- \overline{\nabla}_{e_{l}}k(X,e_{l})\bigg)\\
	&=& -\bigg(  (\delta_{g}k +\text{dtr}_{g}k) \otimes \alpha    \bigg)(X,e_{s}).
\end{eqnarray*}
Thereby we can conclude that  $I=-(\delta_{g}k +\text{dtr}_{g}k) \otimes \alpha$. We compute now $II(X,e_{s})$.
\begin{eqnarray*}
	II(X,e_{s})&=&  \sum_{l\neq s}\text{d}^{\overline{\nabla}}k(X,e_{l},e_{s})\Re e\left\langle  e_{l}\cdot e_{0}\cdot e_{s}\cdot\psi,e_{s}\cdot \psi \right\rangle \\
	&&+ \sum_{m\neq s}\text{d}^{\overline{\nabla}}k(X,e_{s},e_{m})\Re e\left\langle  e_{s}\cdot e_{0}\cdot e_{m}\cdot\psi,e_{s}\cdot \psi \right\rangle \\
	&&+ \sum_{l\neq m,\ l\neq s, \ m\neq s}\text{d}^{\overline{\nabla}}k(X,e_{l},e_{m})\Re e\left\langle  e_{l}\cdot e_{0}\cdot e_{m}\cdot\psi,e_{s}\cdot \psi \right\rangle,\\
\end{eqnarray*}
but the last sum is zero since $\left\langle  e_{k}\cdot e_{0}\cdot e_{m}\cdot\psi,e_{s}\cdot \psi \right\rangle$ is purely imaginary whenever $k,m,s$ are distinct indices. Thereby it comes out that $	II(X,e_{s})=   \big(  \overline{\nabla}_{e_{s}}k(X,\alpha )- \overline{\nabla}_{\alpha}k(X,e_{s} )\big)$, so that we can conclude 
$$II(X,Y)=  \overline{\nabla}_{Y}k(X,\alpha )   -    \overline{\nabla}_{\alpha}k(X,Y ) ,$$
and consequently
\begin{eqnarray*}
V\bigg( \text{Ric}^{g} +2g+ (\text{tr}_{g}k)k - k\circ k \bigg)(X,Y)&=& - \bigg((\delta_{g}k +\text{dtr}_{g}k) \otimes \alpha\bigg)(X,Y)\\
&& + \bigg(  \overline{\nabla}_{Y}k(X,\alpha)   -     \overline{\nabla}_{\alpha}k(X,Y ) \bigg).
\end{eqnarray*}
Moreover the couple $(V,\alpha )$ satisfies the following differential equations
\begin{eqnarray*}
\overline{\nabla}_{X}\alpha (Y)&=&Vk(X,Y) +\frac{{\bf i}}{2}\big((X\cdot Y -Y\cdot X)\cdot \xi ,\xi    \big)\\
\delta^{*}_{g} \alpha &=& Vk \\
\text{d}V(X) &=& k(X,\alpha )+ \textbf{i}\left\langle X\cdot \psi ,\psi \right\rangle 
\end{eqnarray*}
and
\begin{eqnarray*}
(\text{Hess}^{g}V)(X,Y)&=& \overline{\nabla}_{Y}k(X,\alpha)- V(k\circ k)(X,Y)+ Vg(X,Y) + \overline{\nabla}_{X}\alpha (k(Y)) +\overline{\nabla}_{Y}\alpha (k(X))\\
&=&\overline{\nabla}_{Y}k(X,\alpha)-\overline{\nabla}_{\alpha }k(X,Y)- V(k\circ k)(X,Y)+ Vg(X,Y) +\mathcal{L}_{\alpha }k(Y,X)\\
&=&V\bigg( \text{Ric}^{g} +3g+ (\text{tr}_{g}k)k - 2(k\circ k) \bigg)(X,Y)\\
&&+ \bigg((\delta_{g}k +\text{dtr}_{g}k) \otimes \alpha\bigg)(X,Y) + \mathcal{L}_{\alpha }k(Y,X)
\end{eqnarray*}
It is clear that the couple $(V,W):=(V,-\alpha^{\sharp})$ satisfies the first KID equation \cite{BCh}, and using the second KID equation for defining the symmetric tensor $\tau$ of \cite{BCh}, namely
\begin{eqnarray*}
V\left( \tau -\frac{1}{2}(\text{tr}_{g}\tau-\rho ) g \right)&=& V\bigg( \text{Ric}^{g}+ (\text{tr}_{g}k)k - 2(k\circ k) \bigg)-\mathcal{L}_{W}k-(\text{Hess}^{g}V)\\
&=&\bigg((\delta_{g}k +\text{dtr}_{g}k) \otimes W^{\flat }\bigg)-3Vg ,
\end{eqnarray*}
where $2\rho:= \text{Scal}^{g} +(\text{tr}k)^{2}-\left|k\right|^{2}$. Taking the trace of last equation one gets $\text{tr}_{g}\tau-\rho=12$ and consequently
$$V\tau =\bigg((\delta_{g}k +\text{dtr}_{g}k) \otimes W^{\flat }\bigg)+3Vg.$$
Now the equation $\left\langle \widehat{\mathfrak{R}}\xi,\xi\right\rangle =0$ implies $ \bigg( \text{Scal}^{g} +  n(n-1)+ (\text{tr}_{g}k)^{2} - \left|k\right| ^{2}_{g}  \bigg)=2 \left| \delta_{g}k +\text{dtr}_{g}k\right| $. We also know that $ V\bigg( \text{Scal}^{g} +  6+ (\text{tr}_{g}k)^{2} - \left|k\right| ^{2}_{g}  \bigg)= 2\left\langle (\delta_{g}k +\text{dtr}_{g}k),W^{\flat} \right\rangle$, and thereby it is clear that there exists some function on $M$ denoted by $\vartheta  $ such that ${W=\vartheta   ( \delta_{g}k +\text{dtr}_{g}k)}$ and so 
$$V\bigg( \text{Scal}^{g} +  6+ (\text{tr}_{g}k)^{2} - \left|k\right| ^{2}_{g}  \bigg)= 2 \left| \vartheta  \right| \left|\delta_{g}k +\text{dtr}_{g}k\right| ^{2} $$
and therefore 
$$2V  \left| \delta_{g}k +\text{dtr}_{g}k\right| = \bigg( 2\rho  +  6  \bigg) \left| W \right|,$$
which shows that in the Killing development the Killing  vector field $(V,W)$ will be colinear to the cosmological contraints 4-vector $(2\rho +6,2(\delta_{g}k +\text{dtr}_{g}k))$ which is isotropic. It follows that  the Killing vector field $(N,W)$ is also isotropic in the Killing development. We finally obtain the relation $V^{2}(\tau -3g)=\frac{1}{2}(2\rho +6)W^{\flat}\otimes W^{\flat}$ which means that the Killing development is a stationary pp-wave space-time.\\
Supposing furthermore that the constant function $(\xi ,\xi )$ is non-zero, we get ${(\delta_{g} k + \text{dtr}_{g}k)(\xi ,\xi )=0}$ and so $(\delta_{g} k + \text{dtr}_{g}k)=0$  by tracing the equation ${\left\langle \text{d}^{\overline{\nabla}}k(X,Y)\cdot e_{0} \cdot\xi,\xi  \right\rangle=0}$. It follows by the dominant  energy condition that $\text{Scal}^{g} + (\text{tr}k)^{2} - \left|k \right|^{2}=-6$, and we obtain finally that $\tau =3g$.  Thereby the couple $(V,W)$ is a cosmological vacuum KID. It is known \cite{BCh} that in that case $(M,g,k)$ has a cosmological vacuum Killing development denoted by $(\overline{N},\overline{\gamma} )$ which is a stationary 4-dimensional Lorentzian manifold  satisfying $G^{\overline{\gamma}}=3\overline{\gamma}$ and carrying a Killing vector field which is the natural extension of the KID $(V,W)$.\hspace{\stretch{4}}$\square$ \\

\noindent
{\bf Remark.}
It is clear that expecting $m_{0}=0$ so as to define the rigidity situation is much stronger than expecting the degenerate character of $Q$. A good issue would certainly be to use the \textit{geometry at infinity} in the same way as in \cite{ChM1} but in the AdS-asymptotically hyperbolic context, in order to prove under the degenerate character of $Q$ the existence of an isometric embedding of $(M,g,k)$ in AdS.\hspace{\stretch{4}}$\square$\\

\section{Appendix}

\subsection{Non-negativity  of $Q$ seen through its coefficients when $n=3$}

Classical linear algebra results state that every principal minor of $Q$ must be non-negative which give rise to a set of inequalities on the coefficients of $\mathcal{H}$.\\

$$\begin{array}{l}
m_{0}+m_{1} \geq 0 \\
m_{0}-m_{1} \geq 0 \\	
\end{array} $$

$$\begin{array}{rll}
m_{0}^{2}- |m|^{2} &\geq& 0 \\
(m_{0}+m_{1})^{2}-(n_{2}+r_{3})^{2}-(r_{2}-n_{3})^{2} &\geq &0 \\	
(m_{0}-m_{1})^{2}-(n_{2}-r_{3})^{2}-(r_{2}+n_{3})^{2}& \geq& 0 \\	
m_{0}^{2}- m_{1}^{2}- n_{1}^{2}- r_{1}^{2}&\geq& 0 
\end{array} $$

$$\begin{array}{rl}
(m_{0}+m_{1})(m_{0}^{2}- (|m|^{2} +n_{1}^{2}+ r_{1}^{2})) - (m_{0}-m_{1})((n_{2}+r_{3})^{2}+(n_{3}-r_{2})^{2})& \\
 -2((n_{2}+r_{3})(m_{2}n_{1}+m_{3}r_{1})+(-n_{3}+r_{2})(m_{2}r_{1}-m_{3}n_{1}) &\geq 0 \\
 (m_{0}-m_{1})(m_{0}^{2}- (|m|^{2} +n_{1}^{2}+ r_{1}^{2})) - (m_{0}+m_{1})((n_{2}-r_{3})^{2}+(n_{3}+r_{2})^{2})& \\
 +2((n_{2}-r_{3})(m_{2}n_{1}-m_{3}r_{1})+(n_{3}+r_{2})(m_{2}r_{1}+m_{3}n_{1}) &\geq 0 \\
\end{array} $$

$$\begin{array}{ll}
(m^{2}_{0}-(|m|^{2}+|n|^{2}+|r|^{2}))^{2}-4(|m|^{2}|n|^{2}+|m|^{2}|r|^{2}+|n|^{2}|r|^{2})&\\
+4(<m,n>^{2}+<m,r>^{2}+<n,r>^{2})+8m_{0}\text{det}_{\mathbb{R}^{3}}(m,n,r)&\geq 0.\\
\end{array}$$

\subsection{Rigidity Results for the Trautman-Bondi Mass}

Oppositely to the rest of the paper, we consider here the situation of the {\it Trautman-Bondi mass} \cite{ChJL}, namely $(M,g,k)$ is assumed to be Minkowski-asymptotically hyperbolic which means that  the triple $(M,g,k)$ is asymptotic at infinity to a standard hyperbolic slice of Minkowski space-time. It has been proved  (cf. Theorem 5.4 of  \cite{ChJL})  that the {\it Trautman-Bondi four-momentum} $p_{\mu }$ is timelike and future directed under the dominant energy condition (and some other technical assumptions). The aim of this section is to prove some rigidity results for the {\it Trautman-Bondi four-momentum} which are analogous to the statements of section 4. More precisely\\

\noindent
{\bf Theorem.} {\it Under the assumptions of Theorem 5.4 of \cite{ChJL}, and if the component $p_{0}$ of the Trautman-Bondi four-momentum  $p_{\mu }$ vanishes, then $(M,g,k)$ can be isometrically embedded in
Minkowski space-time.}\\

\noindent
{\sc Proof.}
This can be done in the same way as our rigidity theorem: since $p_{\mu }$ is timelike, the condition $p_{0}=0$ implies that $p_{\mu }$ actually vanishes. Consequently there exists a basis of $\nabla $-parallel spinor fields on $M$, where $\nabla $ is the connection on some cylinder $]-\epsilon ,+\epsilon [\times M$ endowed with some Lorentzian metric $\gamma =-\text{d}t^{2}+g_{t}$ (such that $M$ has induced metric $g$ and extrinsic curvature $k$ satisfying the conditions of \cite{ChJL}). Now if $\overline{\nabla }$ denotes the  Levi-Civita connection of $g=g_{0}$, we still have the relation $\nabla _{X}\xi =\overline{\nabla }_{X}\xi - \frac{1}{2}k(X)\cdot e_{0} \cdot \xi$ where $\cdot $ is the Clifford action with respect to $\gamma $. Our spinorial Gauss-Codazzi formula is still valid, that is
$$ R^{\gamma}_{X,Y}= R^{g}_{X,Y} -\frac{1}{2}\left( \text{d}^{\overline{\nabla}}k(X,Y)\cdot e_{0} + \frac{1}{2}\big(k(X)\cdot k(Y)-k(Y)\cdot k(X)\big) \right)\cdot=0 \quad ,$$
and so
\begin{eqnarray*}
 R^{g}&=& \frac{1}{2}k\owedge k \\
\text{d}^{\overline{\nabla}}k&=&0 .
\end{eqnarray*}
Furthermore, the couple $(V,W):=(V,-\alpha^{\sharp})$ is a vacuum KID if one defines $V=<\xi,\xi>$ and the real 1-form $\alpha$  by $\alpha(Y)= \left\langle Y \cdot e_{0} \cdot \xi,\xi \right\rangle$. We consider again the Killing development of $(\widetilde{M},\widetilde{g},\widetilde{k})$ with respect to $(\widetilde{V},\widetilde{W})$, and observe that it must be a geodesically complete stationary solution of vacuum Einstein equations of zero sectional curvature and thereby must be Minkowski space-time (cf. Proposition 23 [p.227] of \cite{O}). Now the same cohomological arguments give $(\widetilde{M},\widetilde{g},\widetilde{k})=(M,g,k)$  which is by construction embedded in its Killing development that is Minkowski. \hspace{\stretch{4}}$\square$\\

\noindent
{\bf Theorem.} {\it Let us suppose that $(M,g,k)$ satisfies the assumptions of Theorem 5.4 of \cite{ChJL} and that $p_{\mu }$ is null. Then there exists some $\nabla$-parallel spinor field $\xi $ such that $\left\langle \mathfrak{R}\xi,\xi\right\rangle =0$ and consequently $(M,g,k)$ is isometrically embeddable in a stationary pp-wave space-time.\\
If furthermore the constant function $(\xi ,\xi )$ is non-zero then $(M,g,k)$ admits a vacuum Killing development which is a stationary solution of the Einstein equations.}\\

\noindent
{\sc Proof.} 
$p_{\mu}$ is null implies the existence of a spinor field $\xi$ satisfying the conditions
\begin{eqnarray*}
	\nabla\xi &=& 0\\
	\left\langle \mathfrak{R}\xi,\xi\right\rangle &=&0 .
\end{eqnarray*}
Then in the same way as in the last Theorem of section 4, but defining here the 2-tensor ${E=:\text{Ric}^{g} +(\text{tr}k)k -k\circ k}$  we obtain that the couple $(V,W)$ is a vacuum KID and the corresponding Killing development satisfies $V^{2}\tau =\rho (W^{\flat}\otimes W^{\flat})$ which means that it is a stationary pp-wave space-time.\\ 
Still using the same computations as in the last Theorem of section 4 and assuming that the constant function $(\xi ,\xi )$ is non-zero we find that the constraints equations are satisfied  (because of the dominant energy condition) and that $\tau =0$. Thereby $(V,W)$ is a vacuum KID and it is known that in this case $(M,g,k)$ has a stationary vacuum Killing development.\\

\hspace{\stretch{4}}$\square$\\

\noindent
{\bf Remark.}
It is clear that expecting $p_{0}=0$ so as to define the rigidity situation is much stronger than expecting the null character of $p_{\mu}$. As in our situation (cf. the remark at the end of section 4), a good issue would certainly be to use the \textit{geometry at infinity} in the same way as in \cite{ChM1} but in the Minkowski-asymptotically hyperbolic context, in order to prove under the equality case of Theorem 5.4 of \cite{ChJL}, the existence of an isometric embedding of $(M,g,k)$ in Minkowski. \hspace{\stretch{4}}$\square$\\

\footnotesize{

{}

\noindent
\textsc{\textbf{Acknowledgements.}} I wish to thank M. Herzlich for his helpful comments as regard the redaction of this text.\\

\noindent
Institut de Mathématiques et de Modélisation de Montpellier (I3M)\\
Université Montpellier II\\
UMR 5149 CNRS\\
Place Eugène Bataillon\\
34095 MONTPELLIER (FRANCE)\\
email:maerten@math.univ-montp2.fr
\\
fax: 33 (0) 4 67 14 35 58\\
tel: 33 (0) 4 67 14 48 47

}


\begin{thebibliography}{99}
\bibitem{And1}
\textsc{M. T. Anderson}, \textit{On stationnary vacuum solutions to the Einstein equations}, Ann. Henri Poincar\'e, \textbf{1} (2000), 5, 977--994.
\bibitem{And2}
\textsc{M. T. Anderson}, \textit{On the structure of solutions to the static vacuum Einstein equations}, Ann. Henri Poincar\'e, \textbf{1} (2000), 6, 995--1042.
\bibitem{AD}
\textsc{L. Andersson, M. Dahl}, \textit{Scalar curvature rigidity for asymptotically locally hyperbolic manifolds}, Ann. Global Anal. Geom., \textbf{16} (1998), 1--27.
\bibitem{B}
\textsc{C. Bär}, \textit{Real spinors and holonomy}, Commun. Math. Phys., \textbf{154} (1993), 509--521.
\bibitem{BGM}
\textsc{C. Bär, P. Gauduchon, A. Moroianu}, \textit{Generalized cylinders in semi-Riemannian and spin geometry}, arXiv:math.DG/0303095 v1.
\bibitem{Bart}
\textsc{R. Bartnik}, \textit{The mass of an asymptotically flat manifold}, Comm. Pure Appl. Math., \textbf{39} (1986), 661--693.
\bibitem{ChB}
\textsc{R. Bartnik, P. T. Chru\'sciel}, \textit{Boundary value problems for Dirac-type equations, with applications}, arXiv:math.DG/0307278v1, 21 Jul 2003.
\bibitem{Baum1}
\textsc{H. Baum}, \textit{Complete Riemannian manifolds with imaginary Killing spinors}, Ann. Global Anal. Geom., \textbf{7} (1989) 4, 205--226.
\bibitem{Baum2}
\textsc{H. Baum}, \textit{Odd-dimensionnal Riemannian manifolds with imaginary Killing spinors}, Ann. Global Anal. Geom., \textbf{7} (1989) 2, 141--154.
\bibitem{BCh} 
\textsc{R. Beig, P. T. Chru\'sciel}, \textit{Killing Initial Data}, Class. Quantum Gravity, \textbf{14} (1997), 1.A, A83-A92.  
\bibitem{Besse}
\textsc{A. L. Besse}, \textit{Einstein manifolds}, Springer.
\bibitem{BT}
\textsc{R. Boot, L. W. Tu}, \textit{Differential Forms in Algebraic Topology}, Graduate Texts in Mathematics, Springer.
\bibitem{BG}
\textsc{J. P. Bourguignon, P. Gauduchon}, {\it Spineurs, opérateurs de Dirac et variations de métriques}, Comm. Math. Phys.,  {\bf 144}  (1992),  n° 3, 581--599.
\bibitem{CGLS}
\textsc{M. Cahen, S. Gutt, L. Lemaire, P. Spindel}, \textit{Killing Spinors}, Bull. Soc. Math. Belg., Ser. \textbf{A 38} (1986), 2, 75--102.
\bibitem{ChH}
\textsc{P. T. Chru\'sciel, M. Herzlich}, \textit{The mass of asymptotically hyperbolic Riemannian manifolds}, Pacific Journal of Mathematics, \textbf{212} (2003), 2, 231--264.
\bibitem{ChJL}
\textsc{P. T. Chru\'sciel, J. Jezierski, S. ~\L\c{e}ski}, \textit{The Trautman-Bondi mass of hyperboloidal initial data sets}, Adv. Theor. Math. Phys., \textbf{8} (2004), 83--139.
\bibitem{ChM}
\textsc{P. T. Chru\'sciel, D. Maerten}, \textit{An upper bound for angular momentum for asymptotically anti-de Sitter space-times}, to appear.
\bibitem{ChM1}
{\sc P. T. Chru\'sciel, D. Maerten }, {\it Killing vectors in asymptotically flat space-times: II. Asymptotically translational Killing vectors and the rigid positive energy theorem in higher dimensions}, to appear.
\bibitem{ChN}
\textsc{P. T. Chru\'sciel, G. Nagy}, \textit{The mass of spacelike hypersurfaces in asymptotically anti-de Sitter space-times}, Adv. Theor. Math. Phys., \textbf{19} (2001), 4, 697--754.
\bibitem{GHW}
\textsc{G. W. Gibbons, C. M. Hull, N. P. Warner}, {\it The stability of gauged supergravity}, Nuclear Phys. B, {\bf  218}  (1983),  n° 1, 173--190. 
\bibitem{HE}
\textsc{S. W. Hawking, G. F. R. Ellis}, \textit{The large scale structure of space-time}, Cambridge University Press.
\bibitem{HT}
\textsc{M. Henneaux, C. Teitelboim}, {\it Asymptotically Anti-de-Sitter Spaces}, Commun. Math. Phys., \textbf{98} (1985), 391--424.
\bibitem{H}
\textsc{M. Herzlich}, \textit{The positive mass theorem for black holes revisited}, J. Geom. Phys., \textbf{26} (1998), 97--111.
\bibitem{LMi}
\textsc{H. B. Lawson, M. L. Michelson}, \textit{Spin Geometry}, Princeton.
\bibitem{L}
\textsc{F. Leitner}, \textit{Imaginary Killing Spinors in Lorentzian Geometry}, J. of Math. Phys., \textbf{44} (2003), no. 10, 4795--4806.
\bibitem{Mo}
\textsc{V. Moncrief}, \textit{Space-time symmetries and linearization stability of the Einstein equations}, J. of Math. Phys., \textbf{16}, 493--498.
\bibitem{O}
\textsc{B. O'Neill}, \textit{Semi-Riemannian Geometry, with applications to Relativity}, Academic Press.
\bibitem{PT}
\textsc{T. Parker, C. Taubes}, \textit{On Witten's proof of the positive energy theorem}, Commun. Math. Phys., \textbf{84} (1982), 223--238.
\bibitem{SY1}
\textsc{R. Schoen, S.-T. Yau}, \textit{On the proof of the positive mass conjecture in general relativity}, Commun. Math. Phys., \textbf{65} (1979), 45--76.
\bibitem{SY2}
\textsc{R. Schoen, S.-T. Yau}, \textit{The energy and the linear momentum of space-times in general relativity}, Commun. Math. Phys., \textbf{79} (1981), 47--51.
\bibitem{SY3}
\textsc{R. Schoen, S.-T. Yau}, \textit{Proof of the positive mass theorem II}, Commun. Math. Phys., \textbf{79} (1981), 231--260.
\bibitem{S}
\textsc{S. T. C. Siklos}, \textit{Lobatchevsky Plane Gravitational Waves}, S.T.C. Siklos, in Galaxies, axisymmetric systems and relativity ed. M.A.H. MacCallum, Cambridge University Press, 1985.
\bibitem{Wald}
\textsc{R. Wald}, \textit{General relativity}, University press of Chicago.
\bibitem{Wang}
\textsc{X. Wang}, \textit{Mass for asymptotically hyperbolic manifolds}, J. Differential Geom., \textbf{57} (2001), 273--299.
\bibitem{Wit}
\textsc{E. Witten}, \textit{A simple proof of the positive energy theorem}, Commun. Math. Phys., \textbf{80} (1981), 381--402.
\bibitem{Z1}
\textsc{X. Zhang}, \textit{Angular momentum and positive mass theorem}, Commun. Math. Phys., \textbf{206} (1999), 137--155.
\bibitem{Z2}
\textsc{X. Zhang}, \textit{A definition of total energy-momenta and the positive mass theorem on asymptotically hyperbolic manifolds}, Commun. Math. Phys., \textbf{249} (2004), 529--548.\\

\end{thebibliography}
\end{document}